\documentclass[12pt,a4paper]{amsart}
\usepackage{amsmath,amsfonts,latexsym,amssymb,graphicx}
\usepackage{amsmath,amsfonts,latexsym,amssymb}

\newcommand{\Z}{{\mathbb Z}}
\newcommand{\N}{{\mathbb N}}

\title[Skew-product for $\cdots$]{Skew-product for group-valued edge labellings of Bratteli diagrams}
\author{A. El Kacimi \& R. Parthasarathy}
\subjclass[2000]{37B10, 54H20}
\keywords{Cantor minimal system, Bratteli diagram, Substitutional system}

\begin{document}
\maketitle

\begin{abstract}
\noindent {\it Group valued edge labellings $\lambda$ of a Bratteli diagram $B$ give rise to a skew-product Bratteli diagram $B(\lambda)$ on which the
group acts. The quotient by the group action of the associated dynamics can be a nontrivial extension
of the dynamics of $B$.  We exhibit a Bratteli diagram for this quotient and construct a morphism to $B$ with unique path lifting property. This is shown to be an isomorphism for the dynamics if a property `loops lifting to loops'
is satisfied by $B(\lambda )\to B$.}

\end{abstract}

\vskip0.5cm

\noindent {\bf 0. Introduction}
\medskip

\noindent In the theory of Cantor minimal dynamical systems Bratteli diagrams are
usually employed to give models.  To associate dynamics to Bratteli
diagrams, in the literature it is always assumed that the edges of the
Bratteli diagram are linearly ordered and moreover invariably one imposes a
further technical assumption of the Bratteli diagram being ``properly
ordered''.  When these assumptions are satisfied the model is the infinite
path space equipped with the so-called `Vershik transformation'.

In this article we associate (subsec.$1.10$) a  Cantor dynamical system to a non-properly
ordered Bratteli diagram.  We use them to construct in subsection 2.2 skew-product dynamical
systems associated to finite group valued edge labellings.  Matui [M] also studied 
`skew-product dynamical systems'
, even though he
did not consider edge labellings.  We were inspired to consider edge labellings
in the context of Bratteli diagrams
%\bigskip
%\hrule
%\bigskip
%{\it Mathematics Subject Classification}: 54H20, 37B10.
%\noindent
%\newpage
by studying the above-mentioned article of Matui and also the article [KP] of Kumjian and Pask who introduced
skew-product $C^*$-algebras associated to group valued edge labellings of
graphs. We remark at the outset that the skew-products which Matui
associates to cocycles arise as our skew-products for particular edge labellings. This is
explained in subsection $2.25$.
For Matui, going modulo the group action on the skew-product system always produces the original dynamical
system. The quotient by $G$ of our skew-product systems would in general be nontrivial extensions of the
original dynamical system. Even for an odometer system our skew-product associated to edge labellings arises
from a cocycle for Matui only when the group is severely restricted. This is described in our example $2.11.2$.
Our main result Theorem $2.12$ has two parts. The first part gives a sufficient condition :-  when a certain
condition ``loops lifting to loops" is satisfied for a Bratteli diagram the quotient by $G$ of our skew-product
is the original system. In the second part of Theorem $2.12$ we describe what happens when the ``loops
lifting to loops" property does not hold; we show that these skew-products also arise from an edge labelling
enjoying the property ``loops lifting to loops" of {\it another Bratteli diagram} which admits a morphism into
the original Bratteli diagram with a certain ``unique path lifting" property.
This is achieved by employing a certain ``tripling" construction. Several examples are given.

We briefly recall some basic concepts and definitions which are fundamental
in the theory of  Cantor dynamical systems.

\medskip

\noindent {\bf 1. Preliminaries}

\medskip

A {\it topological dynamical system} is a pair $(X, \varphi)$ where $X$
is a compact metric space and $\varphi$ is a homeomorphism in $X$.  We say
that $\varphi$ is {\it minimal} if for any $x \in X$, the $\varphi$-orbit
of $x:= \{ \varphi^n (x)\mid  n \in \Z \}$ is dense in $X$.  We say that
$(X, \varphi)$ is a {\it Cantor dynamical system} if
$X$ is a Cantor set, i.e. $X$ is totally disconnected without isolated points.
$(X, \varphi)$ is a {\it Cantor minimal dynamical system} if, in addition,  $\varphi$ is minimal.
Some of the basic concepts of the theory are recalled below, mostly
from the more detailed sources [DHS] and [HPS].

\paragraph*{\bf 1.1. Definition.}  A Bratteli {\it diagram} is an infinite directed
graph $(V,E)$, where $V$ is the vertex set and $E$ is the edge set.  Both
$V$ and $E$ are partitioned into non-empty disjoint finite sets
$$V= V_0 \cup V_1 \cup V_2 \cdots \ {\rm and} \ E=E_1 \cup  E_2 \cup  \cdots $$
There are two maps $r,s:E \to V$ the {\it range} and {\it source} maps.  The following
properties hold:

\begin{enumerate}
\item[(i)] $V_0 = \{v_0\}$ consists of a single point, referred to as the
`top vertex' of the Bratteli diagram
\item[(ii)] $r (E_n) \subseteq V_n, s (E_n) \subseteq V_{n-1}, n=1,2, \cdots
$.  Also $s^{-1} (v) \neq \phi$ $\forall v \in V$ and $r^{-1} (v) \neq
\phi$  for all $v \in V_1, V_2,\cdots$.
\end{enumerate}

Maps between Bratteli diagrams are assumed to preserve gradings and
intertwine the range and source maps.
If $v \in V_n$ and $w \in V_m$, where $m >n$, then a path from $v$ to
$w$ is a sequence of edges $(e_{n+1}, \cdots, e_m)$ such that $s (e_{n+1})
=v, r (e_m)= w$ and $s (e_{j+1})=r (e_j)$. Infinite paths from $v_0 \in V_0$
are defined similarly.  The Bratteli diagram is called {\it simple} if for any $n=0,1,2, \cdots 
$ there exists $m>n$ such that every vertex of $V_n$ can be joined to every vertex
of $V_m$ by a path.

\paragraph*{\bf 1.2. Definition.}  An {\it ordered} Bratteli diagram $(V, E, \geq)$ is a Bratteli 
diagram $(V,E)$ together with a linear order on $r^{-1} (v), \forall v \in
V - \{v_0\} = V_1 \cup  V_2 \cup V_3\cdots $.
We say that an edge $e \in E_n$ is a {\it maximal} edge (resp. {\it minimal} edge)  if $e$ is
maximal (resp. minimal) with respect to the linear order in $r^{-1} (r (e))$.

Given $v \in V_n$, it is easy to see that there exists a unique path
$(e_1, e_2, \cdots ,e_n)$ from $v_0 $ to $v$ such that each $e_i$ is
maximal (resp. minimal).

Note that if $m >n$, then for any $w \in V_m$, the set of paths starting from
$V_n$ and ending at $w$ obtains an induced (lexicographic) linear order:
$$(e_{n+1}, e_{n+2}, \cdots, e_m) > (f_{n+1}, f_{n+2}, \cdots, f_m)$$
if for some $i$ with $n+1 \leq i \leq m, e_j=f_j$ for $ 1 <j \leq m$ and $e_i >f_i$.

\paragraph*{\bf 1.3. Definition.}  A {\it properly ordered} Bratteli diagram is a simple ordered Bratteli diagram
$(V, E \geq )$ which possesses a unique infinite path $x_{\rm max} = (e_1,
e_2, \cdots)$ such that each $e_i$ is a maximal edge and a unique
infinite path $x_{\rm min} = (f_1, f_2, \cdots)$ such that each $f_i$ is a minimal edge.

Given a properly ordered Bratteli diagram $B=(V,E, \geq)$ we denote by $X_B$ its
infinite path space.  So
$$X_B = \{(e_1, e_2, \cdots) \mid e_i \in E_i, r (e_i) = s (e_{i+1}), i= 1,2,
\cdots \}$$
For an initial segment $(e_1, e_2, \cdots,e_n)$ we define the cylinder sets
$$U(e_1, e_2, \cdots e_n) =\{ (f_1, f_2, \cdots ) \in X_B
\mid f_i = e_i, 1 \leq i \leq n \} .$$  By taking cylinder sets to be a basis
for open sets $X_B$ becomes a topological space.  We exclude trivial cases
(where $X_B$ is finite, or has isolated points).  Thus, $X_B$ is a Cantor set.
$X_B$ is a
metric space, where for two paths $x,y$ whose initial segments to level $m$
agree but not to level $m+1, d (x,y)= 1/m+1$.

\paragraph*{\bf 1.4. Definition.} (Vershik map for a properly ordered Bratteli diagram).  If $x= (e_1, e_2, \cdots e_n,
\cdots) \in X_B$ and if at least one $e_i$ is not maximal define
$$V_B (x) =y
=(f_1, f_2, \cdots, f_j, e_{j+1}, e_{j+2}, \cdots ) \in X_B$$
where $e_1, e_2, \cdots , e_{j-1}$ are maximal, $e_j$ is not maximal and has $f_j$ as successor
in the linearly ordered set $r^{-1} (r (e_j))$ and $(f_1,
f_2, \cdots, f_{j-1})$ is the minimal path from $v_0$ to $s (f_j)$.
Extend the above $V_B$ to all of $X_B$ by setting $V_B (x_{\rm max}) = x_{\rm min}$.  Then $(X_B, V_B)$ is a
Cantor minimal dynamical system.

 Next, we describe the construction of a
dynamical system associated to a non-properly ordered Bratteli diagram.
The Bratteli diagram need not be simple.
To motivate this construction, it is perhaps worthwhile to begin by
indicating how it works in the case of an ordered Bratteli diagram
associated to a nested sequence of Kakutani-Rohlin partitions of a Cantor dynamical
system $(X,T)$.

\paragraph*{\bf 1.5. Definition.}  A {\it Kakutani-Rohlin  partition} of the
Cantor minimal system $(X,T)$ is a clopen partition ${\mathcal P}$ of the kind
$${\mathcal P}= \{T^j Z_k \mid k \in A \ {\rm and} \ 0 \leq j <h_k\}$$
where $A$ is a finite set and $h_k$ is a positive integer.  The $k^{th}$
{\it tower} ${\mathcal S}_k$ of ${\mathcal P}$ is $\{T^j Z_k \mid  0 \leq j <h_k\}$ ; its {\it floors} are
$T^j Z_k, (0 \leq j <h_k)$.
The {\it base} of
${\mathcal P}$ is the set $Z=\bigcup_{k \in A} Z_k$.

Let  $\{{\mathcal P}_n\}, (n \in \N)$ be a sequence of Kakutani-Rohlin partitions 
$${\mathcal P}_n =\{ T^j Z_{n,k} \mid k \in A_n,\hbox{and }
0 \leq j < h_{n,k} \} ,$$
with ${\mathcal P}_0 =\{X\}$
and with base $Z_n = \bigcup_{k \in A_n} Z_{n,k}$.  We say that this sequence is
{\it nested} if, for each $n$,
\begin{enumerate}
\item[(i)] $Z_{n+1} \subseteq Z_n$
\item[(ii)]$ {\mathcal P}_{n+1}$ refines the partition ${\mathcal P}_n$.
\end{enumerate}
For the Bratteli-Vershik system $(X_B, V_B)$ of sections 1.3-1.4, one obtains a Kakutani-Rohlin
partition ${\mathcal P}_n$ for each $n$ by taking the sets in the partition to be the
cylinder sets $U(e_1, e_2, \cdots e_n)$ of section 1.3 and taking as the base of the
partition the union $\bigcup U(e_1, e_2, \cdots e_n)$ over minimal paths (i.e.,
each $e_i$ is a minimal edge). This is a nested sequence.

\noindent {\bf 1.6.} To any nested sequence
$\{{\mathcal P}_n\}, (n \in \N)$ of Kakutani-Rohlin partitions we associate an ordered Bratteli diagram $B=(V,E, \geq)$
as follows (see [DHS, section 2.3]): the $\mid A_n \mid$ towers in ${\mathcal P}_n$ are in $1-1$ correspondence
with $V_n$, the set of vertices at level $n$.  Let
$v_{n,k} \in V_n$ correspond to the tower ${\mathcal S}_{n,k} =\{T^j Z_{n,k}
\mid0 \leq j < h_{n,k}\}$ in ${\mathcal P}_n$.  We refer to $T^j Z_{n,k}, 0 \leq j <
h_{n,k}$ as {\it floors} of the tower ${\mathcal S}_{n,k}$ and to $h_{n,k}$ as the {\it
height} of the tower.  {\it We will exclude nested sequences of K-R partitions where the infimum (over $k$ for fixed $n$) of
the height $h_{n,k}$ does not go to infinity with $n$.}
Let us view the tower ${\mathcal S}_{n,k}$ against the partition 
${\mathcal P}_{n-1} =\{T^j Z_{n-1,k} \mid k \in A_{n-1}$, and $0 \leq j < h_{n-1,k}\}$.
As the floors of ${\mathcal S}_{n,k}$ rise from level $j=0$ to level $j=h_{n,k}-1$,  ${\mathcal S}_{n,k}$ will
start traversing a tower ${\mathcal S}_{n-1, i_1}$ from the bottom to the top floor,
then another tower ${\mathcal S}_{n-1, i_2}$ from the bottom to the top floor, then
another tower ${\mathcal S}_{n-1,i_3}$ likewise and so on till a final segment
of ${\mathcal S}_{n,k}$ traverses a tower ${\mathcal S}_{n-1, i_m}$ from the bottom to the top.  Note
that in this final step the top floor $T^j Z_{n,k}$ for $j=h_{n,k}-1$
of ${\mathcal S}_{n,k}$ reaches the top floor $T^q Z_{n-1,i_m}$ for $q= h_{n-1, i_m} -1
$ of ${\mathcal S}_{n-1, i_m}$ as a consequence of the assumption $Z_n\subset
Z_{n-1}$ and the fact that $T^{-1}$ (union of bottom floors) = union of top
floors.  Bearing in mind this order in which ${\mathcal S}_{n,k}$ traverses
${\mathcal S}_{n-1, i_1}, {\mathcal S}_{n-1, i_2}, \cdots, {\mathcal S}_{n-1, i_m}$ we associate $m$ edges,
ordered as $e_{1,k} < e_{2,k} < \cdots < e_{m,k}$ and we set the range and
source maps for edges by $r (e_{j,k})= v_{n,k}$ and $s (e_{j,k}) =v_{n-1, i_j}$.
Note that $m$ depends on the index $k \in A_n $ (and that by convention the
indexing sets $A_n$  are disjoint).  $E_n$ is the disjoint union over
$k \in A_n$ of the edges having range  in $V_n$.

\noindent {\bf 1.7.} For $x \in X$, we define $x_n \in {\mathcal P}_n^{\Z}, n \in \N$ as follows:  $x_n=
(x_{n,i})_{i \in \Z}$, where $x_{n,i} \in {\mathcal P}_n$ is the unique floor in
${\mathcal P}_n$ to which $T^i (x)$ belongs.  If $m >n$, let $j_{m,n} : {\mathcal P}_m \to {\mathcal P}_n$
be the unique map defined by $j_{m,n} (F) =F' $ if $F \subseteq F'$.  (By
abuse of notation, we use the same symbol  $F$ to denote a point of the finite
set ${\mathcal P}_m$ and also to denote the subset of $X$, in the partition ${\mathcal P}_n$,
which $F$ represents).  An important property of the map
$$X \to \prod_n
({\mathcal P}_n^{\Z}), x \mapsto (x_1, x_2, \cdots), x_n = (x_{n,i})_{i \in \Z},$$
defined above is the following:

\noindent {\bf 1.8.} If $F$ and $TF$ are two successive
floors of a ${\mathcal P}_n$-tower and if $x_{n,i}=F$ then $x_{n,i+1}=TF$. If $x_{n,i}$
is the top floor of a ${\mathcal P}_n$-tower, then $x_{n,i+1}$ is the bottom floor of a ${\mathcal P}_n$-tower.
More importantly, given
integers $K$ and $n$, there exist $m>n$ and
a single tower ${\mathcal S}_{m,k}$ of level $m$ such that the finite sequence $(x_{n,i}
)_{-K \leq i\leq  K}$ is  an interval segment contained in
$$\{j_{m,n}(T^\ell
(Z_{m,k})) \mid 0 \leq  \ell < h_{m,k} \}.$$This is a consequence of the assumption that 
the infimum of the heights of level-$n$ towers goes to infinity. It is true that $x_{n,i} =j_{m,n}
(x_{m,i})$, but the sequence $(x_{m,i})_{-K \leq i \leq K}$ need not be an
interval segment of $\{T^{\ell} (Z_{m,k}) \mid 0 \leq \ell <h_{m,k}\}$.

The foregoing observations in the case of an ordered Bratteli diagram associated to a nested sequence
of Kakutani-Rohlin partitions gives us the hint to define a dynamical system  $(X_B,T_B)$ of a non properly ordered
Bratteli diagram $B = (V, E, \geq )$ as follows:

\paragraph*{\bf 1.9. Definition.} For each $n$ define ${\varpi}_n =$ the set of paths from $V_0 $ to $V_n$.  There is an obvious
truncation map $j_{m,n}:{\varpi}_m \to {\varpi}_n$ which truncates paths from $V_0$
to $V_m$ to the initial segment ending in $V_n$. For
each $v \in V_n$, the set $\varpi (v)$ of paths from $\{ \ast \} \in V_0$ ending at $v$ will be called a
`${\varpi}_n$-{\it tower parametrised by} $v$'; the cardinality $|\varpi (v)|$ will be referred to as
the height of this tower.  Each tower is a linearly ordered
set (whose elements may be referred to as floors of the tower) since
paths from $v_0$ to $v$ acquire a linear order ({\it cf.} 1.2). {\it We will exclude unusual examples
of ordered Bratteli diagram where the infimum of the height of level-$n$ towers does not go to infinity with $n$},
(for example like [HPS, Example 3.2]).
Now, we define

\noindent{\bf 1.10.  Definition.} $X_B = \{x = (x_1, x_2, \cdots, x_n, \cdots ) \}$ where
\begin{enumerate}
\item[(i)] $x_n = (x_{n,i})_{i \in \Z} \in {\varpi}_n^{\Z}$,
\item[(ii)]  $j_{m,n} (x_{m,i}) =x_{n,i}$ for $m>n$ and $i \in \Z$    and
\item[(iii)] given $n$ and
$K$ there exists $m$ such that  $m >n$ and a vertex $v \in V_m$,
such that the interval segment $x_n [-K, K]:= (x_{n,-K},
x_{n,-K+1}, \cdots, x_{n,K})$ is obtained by applying $j_{m,n}$ to
an interval segment of the linearly ordered set of paths from
$v_0$ to $v$.
\end{enumerate}
The condition (iii) is the crucial part of the definition. Without it what one 
gets is an inverse system.

The condition (iii) implies that a property similar to (1.8) holds.
Since each ${\varpi}_n$ is a finite set ${\varpi}_n^{\Z}$ has a product topology which
makes it a compact set - in fact a Cantor set.  Likewise, $\prod_n ({\varpi}_n^{\Z})$
is again a Cantor set.  Thus, $X_B \subseteq \prod_n ({\varpi}_n^{\Z})$ has
an induced topology. The lemma below and the following proposition are 
analogous to corresponding facts for the Vershik model associated to properly 
ordered Bratteli diagrams.

\smallskip
\paragraph*{\bf 1.11. Lemma.}  {\it The topological space $X_B$ is compact.}
\smallskip
\paragraph{\bf{Proof.}} We show that $X_B$ is closed in $\prod_n {\varpi}_n^{\Z}$.  Let $z=(z_1,
z_2, \cdots, z_n, \cdots )$, where $z_n \in {\varpi}_n^{\Z}$.  Assume that
$z=\lim_{m \to \infty} w_m$, where $w_m \in \prod_n {\varpi}_n^{\Z}
$ and moreover $w_m \in X_B$. Let $K_1, K_2, \cdots,$ be a strictly increasing
sequence of positive integers.  Define neighbourhoods $U_1, U_2, \cdots,
U_m, \cdots$ of  $z$ in $\prod_n {\varpi}_n^{\Z}$ shrinking to $z$ by
$U_m = \left\{ z' = (z'_1, z'_2, \cdots, z'_m, \cdots ) \in \prod_n {\varpi}_n^{\Z}
\right\} $
where $z'_k$ has the same coordinates as $z_k$ in the range $[-K_m, K_m ]$, i.e.
$z'_k [-K_m, K_m]= z_k [-K_m, K_m]$ for $1 \leq k \leq m$.  Since $z=
\lim w_j$, for any given $m$, $\exists J(m)$, such that $``j>J(m)'' \Rightarrow
w_j \in U_m$.  But, since $w_j \in X_B$, we conclude that for $1 \leq k \leq
m$, the interval segments $z_k [-K_m, K_m]$ are obtained by applying
$j_{M,k}$ to an interval segment of the sequence of floors of a single
${\varpi}_M-$tower (for some $M)$.  This shows that $z\in X_B$.  \hfill $\square $

Denote by $T_B$ the restriction of the shift operator to $X_B$. 
So, if $ x = (x_1, x_2, \cdots, x_n, \cdots ), $ where $x_n=(x_{n,i})_{i \in \Z} \in {\varpi}_n^{\Z}$,
then $T_B(x) = (x_1',x_2',\cdots,x_n',\cdots) $,
where $x_n'=  (x_{n,i}')_{i \in \Z} \in {\varpi}_n^{\Z}$ and $x_{n,i}'=x_{n,i+1}$.

$(X_B, T_B)
$ will be called the dynamical system associated to $B=(V, E, \geq)$.

\smallskip
\paragraph*{\bf 1.12. Proposition.} {\it If $B=(V, E, \geq)$ is a simple ordered
Bratteli diagram, then $(X_B, T_B)$ is a Cantor minimal dynamical system.}
\smallskip
\paragraph*{\bf{Proof.}}  Let $x,y \in X_B$, where $x= (x_1,  \cdots,
x_m, \cdots)$ and $y = (y_1, \cdots, y_{m},\\ \cdots )$ satisfy the
conditions i), ii) and iii) of definition (1.10).  We will show that $x$ belongs
to the closure of the orbit $\{T_B^i (y) \mid i \in \Z\}$. For integers $n$ and $K$
define the neighbourhood $U(n,K,x)$ of $x$ to be the set\\
$\{z = (z_1, z_2, \cdots, z_m, \cdots) \in X_B \mid
z_m [-K, K] = x_m [-K, K]$ for $1 \leq m \leq n\}$.\\
Given any neighbourhood
$U$ of $x$, choose $n$ and $K$ such that
$U \supseteq U (n,K,x)$.
Since $x \in X_B$, we
can choose $J$ and a ${\mathcal P}_J$-tower ${\mathcal S}_v$ of paths from $v_0$ ending at
a fixed $v \in V_J$, such that  $x_n [-K, K]$ is obtained by applying
$j_{J,n}$ to an interval segment of ${\mathcal S}_v$.  Since the Bratteli diagram
$(V, E, \geq )$ is simple, we can choose $L >J$ such that every point
of $V_J$ is connected to every point of $V_L$ by a path.   This implies
that $x_n [-K,K]$ occurs as an interval segment of the sequence
obtained by applying $j_{L,n}$ to the sequence of floors of {\it any}
${\mathcal P}_L$-tower.  Now, choose a $<b$ such that $y_L [a,b]$ is  the sequence
of all floors of a ${\mathcal P}_L$-tower.  From the preceding observation, we can choose
$d$ such that $a \leq d < d+2K \leq b$ and $x_n [-K,K] =j_{L,n} (y_L [d,d+
2K])$.  But since $y \in X_B$, $y$ satisfies condition (ii)  in definition 1.10.
So $x_n [-K,K] = y_n[d, d+2K]$ and we conclude $T_B^{d+K}
(y) \in U (n,K,x)$. \hfill $\square $

In (1.7), given a nested sequence of Kakutani-Rohlin partitions of
$(X,T)$, we defined a map from $(X,T)$ to the dynamical
system $(X_B, T_B)$ of the associated ordered Bratteli diagram.  It follows
that if $(X,T)$ is minimal, and if the Bratteli diagram of the nested
sequence of K-R partitions is a simple Bratteli diagram, then
$(X,T) \to (X_B, T_B)$ is onto. If the topology of $(X,T)$ is spanned
by the collection of the clopen sets belonging to the K-R partitions
then clearly the map $(X,T) \to (X_B, T_B)$ is injective. In particular,
if the Bratteli diagram is properly ordered then the Bratteli-Vershik
system is naturally isomorphic to the system given by our construction
in 1.10.

\smallskip

Note that the same term {\it `towers'} has been used to denote two separate
but related objects [in (1.5) and (1.9)].
For $v \in V_n $,
let $y$ be a path from $\{ *\}$ to $v$ in $(V, E,
\geq)$.  So, $y$ is a  `floor' (consisting of the single element $y$)  belonging to the ${\varpi}_n$-
tower $\varpi (v)$ (a finite set) parametrized by $v\in V_n$ - all in the sense of  $(1.9)$.
Here, $\varpi (v)$= all paths from $\{ *\}$ to $v$.
Put ${\mathcal F}_y = \{x= (x_1, x_2, \cdots, x_n, \cdots ) \in X_B  \mid
x_{n,0} =y \}$. ${\mathcal F}_y$ is a clopen set of the Cantor set $X_B$.  Put
${\mathcal P}_n = \{{\mathcal F}_y \mid y \in \varpi (v), v \in V_n\}.$
Then, in the sense of $(1.5)$ ${\mathcal P}_n$ is a  K-R partition of $X_B$ whose 
base is the union of $\bigcup {\mathcal F}_y,( y \text{ minimal } \in \varpi (v), v \in V_n) $.
Its towers
${S}_v$ are parametrized by $v\in V_n$: ${S}_v = \{{\mathcal F}_y\mid  y\in \varpi (v)\}$.
${\mathcal F}_y,(y\in \varpi (v))$ are the floors of the tower ${S}_v$.
(We encountered this K-R partition earlier in the case of the Bratteli-Vershik system at the end of
1.5).It is easy to see that the ordered Bratteli diagram obtained from $\{{\mathcal F}_y \mid y \in \varpi (v), v \in V_n\}$
is $(V,E,\geq)$. (compare the last two lines in the proof of [HPS, Theorem 4.5] and also the opening observation
in the proof of [DHS, Proposition 16]).

\medskip
\noindent {\bf 2. Skew-product dynamical systems}
\medskip

\noindent {\bf 2.1.} We recall the notion of `skew-product dynamical systems' which was defined
by Matui in [M].  Let $(X,T)$ be a Cantor dynamical system.  Let $G$ be
a finite group and $c: X \to G$ a continuous function.  We set $Y= X \times G
$ and define a continuous map $S: Y \to Y$ by $S (x,g) = (T(x), g c (T(x)))$.
The dynamical system $(Y,S)$ is called the {\it skew-product extension
associated} to the $G$-valued {\it cocycle} $c$.  For $g \in G$, define
$\gamma_g :Y \to Y$ by  $\gamma_g (x,h) = (x,gh)$.  Clearly, $\gamma_g \circ
S =S \circ \gamma_g, \forall g$.  Furthermore, the map
$X \times G \to X$ given by projection to the first factor  displays $(Y,S)$
as an extension of $(X,T)$, (or, $(X,T)$ as a factor of $(Y,S))$.
In a different context, Kumjian and Pask defined in [KP] the $C^\ast $-algebra
associated to a graph and further when the edges of a graph are labelled with
values in a finite group $G$, they considered a `skew-product graph' and
`skew-product $C^*$-algebra' associated to the labelling.

Motivated by these two works, we employ our construction of subsection $1.10$
to associate a `skew-product' dynamical system when the edges of an
ordered Bratteli diagram are labelled with values in a finite group
$G$.

\noindent {\bf 2.2.}  Let $(V, E, \geq)$ be an ordered Bratteli
diagram.  Let $G$ be a finite group and $\lambda: E \to G$ a map
defined on the edges with values in $G$
($\lambda{\overset{defn.}=}$ the `{\it labels}' on edges).  We
define a new Bratteli diagram $B(\lambda){\overset{defn.}=}
(V_{\lambda}, E_{\lambda})$ as follows: $V_{0,\lambda} = \{\ast
\}, V_{n,\lambda} = V_n  \times G, (n \geq 1)$ and $E_{n,\lambda}
=E_n \times G, (n \geq 1)$.  The source and range maps on
$E_{\lambda}$ are defined by $r (e,g)= (r(e),g)$, and $s (e,g)=
(s(e), g \lambda (e))$.  If $e \in E_1$, we define $s(e,g)
=\{*\}$. We define a map $\pi : (V_{\lambda}, E_{\lambda}) \to
(V,E)$ by $\pi (v_{\lambda, 0}) =v_0, \pi (v,g)=v$ if $v \in V_n,
n \geq 1$ and $\pi (e,g)= e $.  One sees that $\pi$ commutes  with
the range and source map.  It is easy to see that $\pi  \vert
_{r^{-1} (v,g)}$ is a bijection onto $r^{-1} (v), \forall v \in
V-\{v_0\}$ and $\forall g \in G$.  Thus if $(V, E, \geq)$ is an
ordered Bratteli diagram there is a unique order in $E_{\lambda}$
such that $\pi : (V_{\lambda}, E_{\lambda}, \geq ) \to (V, E, \geq
)$ is order preserving.  In the sequel we assume that
$(V_{\lambda}, E_{\lambda})$ is equipped with this order. The
dynamical system $(X_{\lambda},T_{\lambda}){\overset{defn.}=}$ the
system constructed in $1.10$ for this ordered Bratteli diagram
$B(\lambda)=(V_{\lambda}, E_{\lambda},\geq)$ will be called the
{\it skew-product system} for the edge labelling $\lambda$. We
also remark that
$$\pi: (V_{\lambda}, E_{\lambda}) \to (V,E) \text{ has the
`{\it unique path lifting}' property }
$$
in the following sense.  If $m > n \geq 1$,
and $(e_n, e_{n+1}, \cdots, e_m)$ is a path in $(V,E)$ from $V_{n-1}$
to $V_m$ with $r(e_m)=v$ then for any $g \in G$, there is a unique path
$(\tilde{e}_n, \tilde{e}_{n+1}, \cdots, \tilde{e}_{m})$ in $(V_{\lambda},
E_{\lambda})$ which maps onto $(e_n, e_{n+1}, \cdots, e_m)$
under $\pi$ and such that $r (\tilde{e}_m) = (v,g)$.

\noindent [\underbar{CAUTION} :- Even in the presence of unique path lifting property
two different edges on the left with the same source may map into the same edge on the right.  See example (2.8)
below. What the
property asserts is that two different edges on the left with the same range cannot map to the same
edge on the right.]

The group $G$ acts on
$(V_{\lambda}, E_{\lambda}, \geq)$ by $\gamma_g (v,h)= (v,gh) $ and
$\gamma_g (e,h) = (e,gh)$  for $v \in V-\{v_0\}$ and of course $\gamma_g
(v_{0,\lambda}) = v_{0,\lambda}$.

Given the labelling $\lambda$ we can extend the labelling to paths from
$V_{n-1}$ to $V_m$, for $m >n$.  With notation as above, we define $\lambda
(e_n, e_{n+1}, \cdots, e_m)= \lambda (e_m) \lambda (e_{m-1}) \cdots
\lambda (e_n)$.  Let $ \{n_k\}^{\infty}_{k=0}$ be a subsequence of $\{0,
1,2, \cdots \}$ where we assume $n_0=0$.  A Bratteli diagram $(V',E')$
is called a {\it `telescoping'} of $(V,E)$ if $V'_k=V_{n_k}$ and
$E'_k$ consists of paths $(e_{n_{k-1}+1}, \cdots, e_{n_k})$ from
$V_{n_{k-1}}$ to $V_{n_k}$ in $(V,E)$, the range and source maps being the obvious
ones.  Thus, for every telescoping the labelling $\lambda$ on the
edges of $(V,E)$ gives rise to a labelling on the edges of $(V', E')$.

\noindent {\bf 2.3. Remark.}  Suppose there is a telescoping $(V',E')$ of $(V,E)$
such that the induced labelling $\lambda$ on $(V',E')$ has the following
property: given $k$, there exists $v' \in V'_k$ and $w' \in V'_{k+1}$,
such that $\lambda (s^{-1} (v') \cap r^{-1} (w'))=G$.  If in addition
$(V,E)$ is simple, then  $(V_{\lambda}, E_{\lambda})$ is simple.

\noindent {\bf 2.4. Stationary Bratteli diagrams.}
A Bratteli diagram is stationary if the diagram repeats itself
after level $1$. (One may relax by allowing a period from
some level onwards; but, a telescoping will be stationary
in the above restricted sense.) If $(V,E,\geq )$ is an ordered Bratteli diagram and the
diagram together with the order repeats itself after level $1$, then
$(V, E, \geq)$ will be called a stationary ordered Bratteli diagram.
We refer the reader to [DHS, section (3.3)] for the usual definition of
a substitutional system and how they give rise to stationary Bratteli diagrams.
Some details are recalled below.
Let $(V,E, \geq)$ be as above and suppose moreover that it is a simple Bratteli
diagram and that $\lambda$ is a labelling of the edges with values
in a finite group $G$.  Assume that  the labelling is stationary: so, we have

\begin{enumerate}
\item[ (1)] an enumeration $\{v_{n,1}, v_{n,2}, \cdots, v_{n,L}\}$ of $V_n, \forall
n \geq 1$,

\item[ (2)]  for $n>1$ and $1 \leq j \leq L$ an enumeration
$\{e_{n,j,1}, e_{n,j,2}, \cdots, e_{n,j,a_j}\} $ of $\ r^{-1} (e_{n,j})$
which is assumed to be listed in the linear order in $r^{-1} (v_{n,j})$,

\item[ (3)]  in the enumerations above, $L$ does not depend on $n$ and $a_j$ depends
only on $j$  and not on $n$.   Moreover, if $n, m >1$, if $ 1 \leq j \leq L,
1 \leq  k \leq L, 1 \leq i \leq a_j$, then $`` s (e_{n,j,i}) = v_{n-1,k}$'' $\Longrightarrow ``s (e_{m,j,i}) = v_{m-1, k}$'',

\item[ (4)] with notation as in 3) above, $\lambda (e_{n,j,i}) = \lambda (e_{m,j,i})$.
\end{enumerate}

If $S$ is a set of generators for $G$ and if in addition to the above,
we also assume that
$\lambda (\{r^{-1} (v_{n,j})\} \cap \{s^{-1} (v_{n-1,k})\}) \supseteq S \cup \{e\},
\forall j,k$ between 1 and $L$ then $ (V_{\lambda}, E_{\lambda})$ is
simple.

\noindent{\bf 2.5. Substitutional systems.} Let $A$ be an alphabet set.
Write $A^+$ for the set of words of finite length in the alphabets of $A$.
Let $\sigma:A \to A^+ $ be a substitution, written,
$\sigma (a) = \alpha \beta \gamma \cdots $
The stationary ordered Bratteli diagram $B=(V,E,\geq )$ associated to $(A, \sigma)$ ({\it cf.} [DHS, section 3.3] can be
described as
$$V_n =A , \forall \text{ }n \geq 1, V_0 =\{*\}$$
$$E_n =\{(a,k,b) \mid  a,b \in A, k \in \N,\text{ } a \text{ is the }k^{th} \text{ alphabet
in the word }\sigma  (b)\}.$$
(The reader who prefers a more carefully evolved notation can consider intoducing
an extra factor `$\times \{n\}$' so that vertices and edges at different levels are
seen to be disjoint).
The source and range maps $s$ and $r$ are
defined by $s (a,k,b)=a, r (a, k,b)=b$.  In the linear order in $r^{-1}
(b)$, $ (a,k,b)$ is the $k^{th}$ edge.  If $\lambda$ is a stationary labelling
on edges with values in a finite group $G$, then define a new `skew-product' substitutional
system $(A_{\lambda}, \sigma_{\lambda})$ as follows:=
\begin{enumerate}
\item[]$A_{\lambda} =A \times G$
\item[]$E_{n,\lambda} =\{[(a,g),k, (b,h)] \mid \text{ (i) } (a,k,b) \in E_n,\text{ (ii) } g=h \cdot  \lambda (a,k,b) \}$.
\end{enumerate}
Define
$\sigma_{\lambda}: A_{\lambda} \to A^+_{\lambda}$ by the rule that
$$(a,g) \text{ is the } k^{{\it th}} \text{ alphabet in the word  }\sigma_{\lambda} (b,h)$$
if
$$ \ [(a,g),k, (b,h)] \in E_{n,\lambda}.$$
Set $V_{\lambda}=\{*\}\cup \{\bigcup_{n\in \N} V_{n,\lambda}\}$ where each $V_{n,\lambda}$ is a (disjoint) copy of $A_{\lambda}$
and likewise, $E_{\lambda}= \cup_n E_{n,\lambda}$. The source and range maps are defined by
$s([(a,g),k, (b,h)])=(a,g)$ and $r([(a,g),k, (b,h)])=(b,h)$.
Then $(V_{\lambda}, E_{\lambda}, \geq )$ is the stationary ordered
Bratteli diagram arising from the substitutional system $(A_{\lambda},
\sigma_{\lambda})$. The group $G$ acts on $(V_{\lambda}, E_{\lambda}, \geq )$ by
$$g(a,g_1) = (a, gg_1)$$
$$g[(a,g_1), k, (b,g_2)]= [(a,gg_1), k, (b,gg_2)].$$
The ordered Bratteli diagram $(V_{\lambda},E_{\lambda},\geq )$ thus obtained is the same as the skew-product of (2.2)
if one starts with the $(V,E,\geq)$ in the beginning of (2.5).

\smallskip

To the stationary ordered Bratteli diagram $B$ of $(A,\sigma)$ (which may not be properly ordered unless
$\sigma$ is a primitive, aperiodic, proper substitution, -- see [DHS, section 3]) we can associate
a dynamical system $X_B$ following the construction of 1.10; we will see that this is naturally
isomorphic to the substitutional dynamical system $X_{\sigma}$ associated to $(A,\sigma)$ defined for example
in [DHS, section 3.3.1]. Let us use the nested sequence $({\mathcal P}_n)$ of Kakutani-Rohlin partitions in $X_{\sigma}$
defined in [DHS, corollary 13]; the ordered Bratteli diagram associated by section 1.6 to this nested 
sequence $({\mathcal P}_n)$ is the same as the stationary ordered Bratteli diagram $B$ in the beginning
of this section. (this is the opening observation in the proof of [DHS, Proposition 16], but, this observation
also remains valid for improper substitutions).
Next apply the map $X_{\sigma} \rightarrow X_B$ that was defined in 1.7. As 
remarked after the proof of Proposition 1.12, this map is onto.  That this map is also injective can
be seen as follows: let $x\in X_{\sigma}$. So, $x=(a_i)_{i\in \Z} \in A^{\Z}$. Let the image of $x$
in $X_B$ by the map of 1.7 be denoted by $(x_1,x_2,\cdots )$ where $x_n = (x_{n,i})_{i \in \Z} \in
{\mathcal P}_n^{\Z}$. The sequence $x = (a_i)_i  \in A^{\Z}$ can be read off from the first co-ordinate $x_1
=( x_{1,i})_i \in {\mathcal P}_1^{\Z} $. 
Choose any element $z\in x_{1,i}$.
Writing, uniquely, as in [DHS, Corollary 12(ii)] $z = T_{\sigma}^k(\sigma (y))$ 
where $y\in X_{\sigma}$ and $0\leq k < |\sigma(y_0)|, $ it is at once seen that
$a_i  \text{ is the }k^{th} \text{
alphabet in } \sigma (y_0)$ and is independent of $z \in x_{1,i}$. In fact
it is a part of the definition [DHS, Corollary 13] of the Kakutani-Rohlin partition ${\mathcal P}_n$ that
$y_0 \text{  and }k$ are independent of $z\in x_{1,i}$

\noindent{\bf 2.6. Example.} Consider the (proper) substitutional (dynamical)
system in the alphabet set $A= \{ X,Y\}$ given by $\sigma (X) =X X Y$,
$\sigma (Y) = XYY$.  The corresponding (stationary) Bratteli
diagram $(V, E, \geq)$ can be represented by
%\newpage
%\medskip

\centerline{\bf Diagram 1}
\begin{center}
\includegraphics[totalheight=10cm]{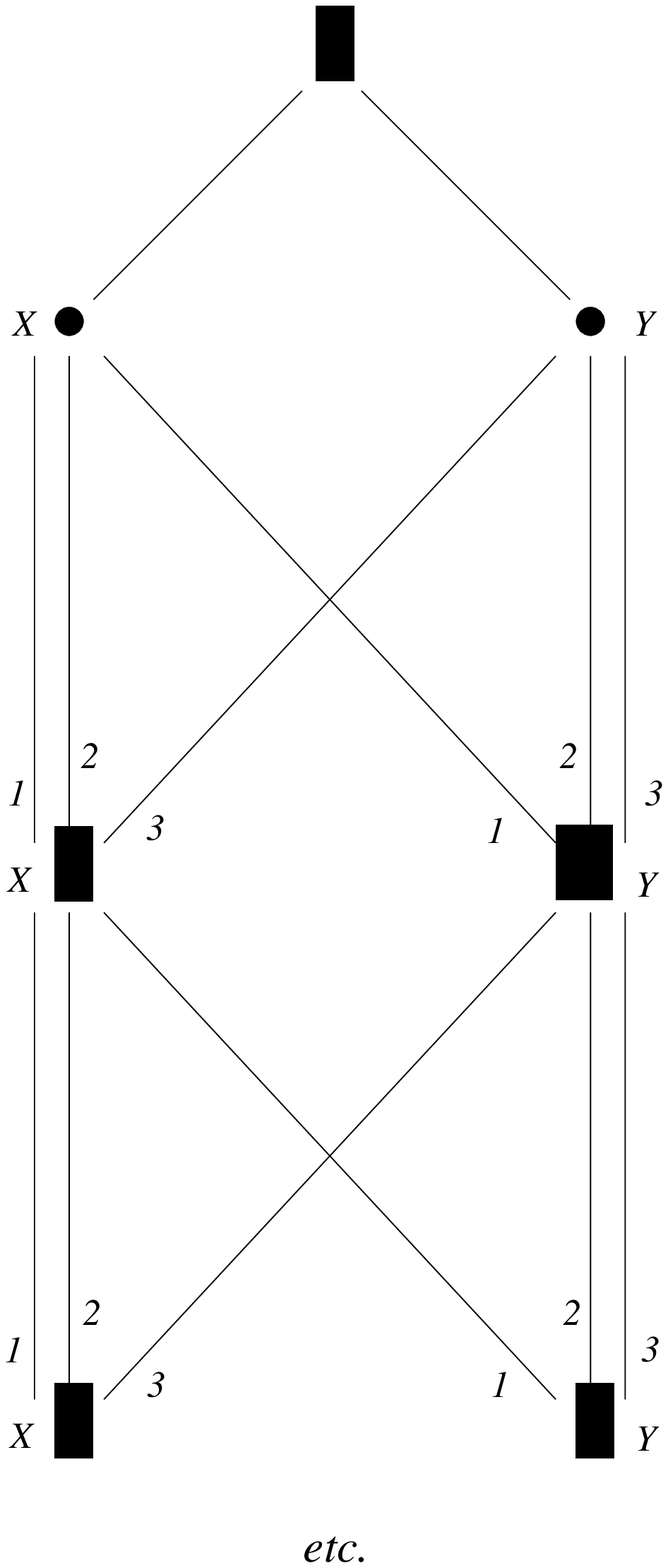}
\end{center}

%\vskip 7cm

\noindent{\bf 2.7. Example $(2.6)$ continued:}
For our substitutional system $A= \{X,Y\}$, $\sigma: A \to A^+$ given by
$\sigma (X)=XXY, \sigma (Y) = XYY$, the associated Bratteli diagram is
stationary.  Each $V_n$ can be identified with $\{X,Y\}$.  There are six
edges connecting $V_n$ and $V_{n+1}$ noted $e^1_{XX}, e^2_{YX},
e^3_{YX},\\ e^1_{XY}, e^1_{XY}, e^2_{YY}, e^3_{YY}$ with source and
range maps defined by
$$
\begin{array}{lll}
s(e^1_{XX})=X &,& r (e^1_{XX})=X \\
s(e^2_{XX}) =X &,& r (e^2_{XX}) =X \\
s(e^3_{YX}) =Y &,& r (e^3_{YX}) =X \\
s(e^1_{XY}) =X &,& r (e^1_{XY}) =Y \\
s(e^2_{YY}) =Y &,& r (e^2_{YY}) =Y \\
s(e^3_{YY}) =Y &,& r (e^{3}_{YY}) =Y
\end{array}
$$

The edges with range $X$ are linearly ordered as $\{e^1_{XX}, e^{2}_{XX},
e^{3}_{YX}\}$.  The edges with range $Y$ are linearly ordered as
$\{e^1_{XY}, e^{2}_{YY}, e^{3}_{YY} \}$.  Consider the stationary labelling
with values in the group $G=\Z_2 = \{ \overline{0}, \overline{1}\}$
given by $\lambda (e^{2}_{YY})=
\overline{1}$ and $\lambda (f)= \overline{0}$ for the remaining
edges.

We now describe $(V_{\lambda}, E_{\lambda})$.  Recall $V_{n,\lambda}=
V_n \times G$. Thus each $V_{n, \lambda}$ can be identified with
$\{X,Y\} \times  \Z_2$.  We identify $\{X,Y\}  \times \Z_2$ with the set of alphabets
$\{x,x',y, y'\}$ by $x= (X, \overline{0}), x'= (X, \overline{1}), y= (Y,
\overline{0}), y' = (Y, \overline{1})$.  Then the edges $V_{n, \lambda} -
V_{n+1, \lambda}$ can be diagrammatically represented by
\medskip
\centerline{\bf Diagram 2}
\begin{center}
\includegraphics[totalheight=8cm]{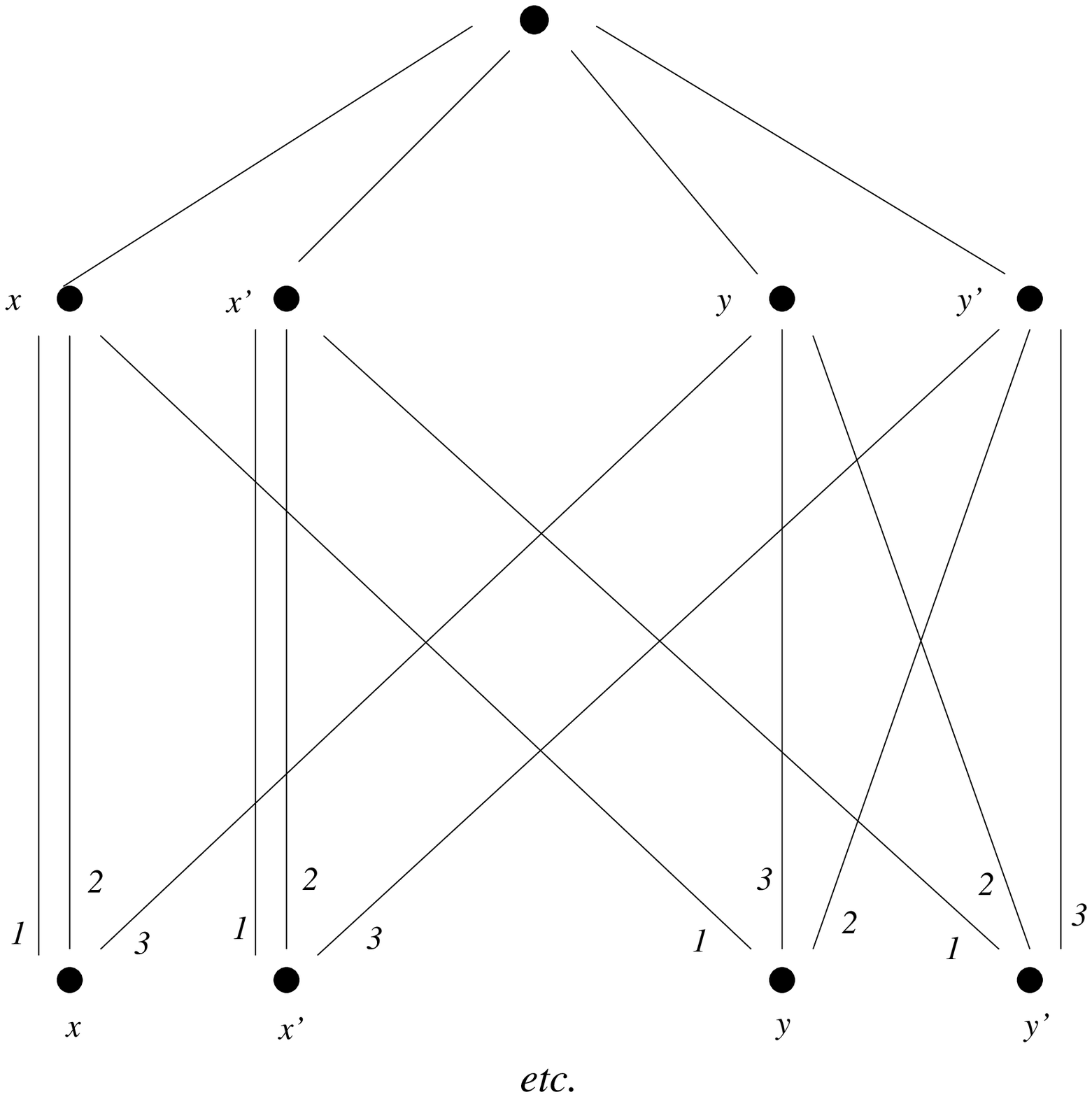} \end{center}
%\vskip6cm

One concludes that  $(V_{\lambda}, E_{\lambda}, \geq)$ is the stationary ordered Bratteli diagram associated to the substitutional system on the alphabet
set $A_{\lambda}=\{x,x',y, y'\}$ where the substitution $\sigma_{\lambda}: A_{\lambda} \rightarrow A_{\lambda}^+ $ is given by
$$
\begin{array}{lll}
x & \mapsto & xxy \\
x' & \mapsto & x'x'y' \\
y & \mapsto & xy'y \\
y' & \mapsto & x'yy'
\end{array}
$$
The Bratteli diagram $(V_{\lambda}, E_{\lambda}, \geq)$ is simple:- any vertex at level $n$ can be connected to any vertex at
level $n+3$ by a path.

\noindent{\bf 2.8. Example where $G$ is infinite. }
Let $G$ be a residually finite group. Let $G  \supset
N_0 \supset N_1 \supset \cdots \supset N_i \supset \cdots$ be a decreasing
sequence of cofinite normal subgroups of $G$ such that $\bigcap_i (N_i) = \{e\}$.
Let $\varphi_i:G\to G/N_i$ be the canonical projection.
We use freely the notation for a substitutional system $(A,\sigma)$ set up in $(2.5)$.
Let $\lambda$ be a stationary labelling on $E = \{(a,k,b)\mid a \text{ is the }k^{th} \text{
alphabet in } \sigma (b)\}$ with values in $G$. The skew-product system in this example is
defined as follows:-
$$ V_{n,\lambda} = A\times (G/N_n), \forall n$$
$$ E_{n,\lambda} = \{[(a,\varphi_{n-1}(g)),k,(b,\varphi_n(h))]\} $$
where
$a,b\in A, (a,k,b)\in E$ and $\varphi_{n-1}(g)=\varphi_{n-1}(h\cdot \lambda (a,k,b))$.

In the statements below, we do not assume that $(V_{\lambda}, E_{\lambda},
\geq)$ arises  from a stationary labelling of a stationary ordered Bratteli
diagram.  Nor do we assume that $(V_{\lambda}, E_{\lambda},
\geq)$ is simple. We  only assume that

\begin{enumerate}
\item[(1)] $(V, E, \geq )$ is a simple ordered Bratteli
diagram and
\item[(2)] $\lambda$ is a labelling with values in $G$.
\end{enumerate}

\noindent{\bf 2.9.} Let $(X,T)$ and $(X_{\lambda}, T_{\lambda})$ be the Cantor dynamical
systems associated respectively to $(V,E,\geq )$ and $(V_{\lambda}, E_{\lambda}, \geq )$ as in subsection (1.10). The map $\pi: (V_{\lambda}, E_{\lambda}, \geq ) \to
(V,E, \geq) $ sends paths    from
$\{ *\}$ to $(v,g)$ in $(V_{\lambda}, E_{\lambda})$ to paths   from $\{*\}$ to $v$ in $(V, E)$ and respects truncation
({\it cf}. 1.9).
The unique path lifting property   implies that $\pi$ maps a  ${\varpi}_{n, \lambda}$-
tower in $(V_{\lambda}, E_{\lambda}, \geq )$ parametrised by $(v,g)
\in V_{n,\lambda}$ bijectively onto the ${\varpi}_n$-tower  in $(V, E, \geq )$
parametrised by $v \in V_n$ (and of course, respects the linear order
of the floors).  Thus $\pi$ induces a map
$\pi: (X_{\lambda}, T_{\lambda}) \to (X,T)$
between the two dynamical systems.  The action of $G$ on $(V_{\lambda},
E_{\lambda}, \geq)$ given by  $\{\gamma_g\}_{g \in G}$ gives rise to a
free action $\{\gamma_g\}_{g \in G}$ of $G$ on $(X_{\lambda}, T_{\lambda})$
such that $\pi \circ \gamma_g = \pi$.

\noindent{\bf 2.10. Remark.}  As observed by Matui [M], when $G$ acts freely on
a Cantor dynamical system $(Y,S)$, the latter is isomorphic to the skew-product
dynamical system arising from a  $G$-valued cocycle on the quotient
of the action.  If one starts with a cocycle on $(X,T)$ then for the
associated skew-product system $(Y,S)$ where $Y=X \times G$, the
quotient is naturally isomorphic to $(X,T)$.  But, in our situation
$(X_{\lambda}, T_{\lambda})$ arises from a $G$-valued labelling of the edges of $(V,E,\geq)
$, the dynamical system associated to $(V,E, \geq)$ is $(X,T)$, but
the quotient of $(X_{\lambda}, T_{\lambda})$ by the $G$-action
$\{\gamma_g\}_{g \in G}$
may still be a nontrivial extension of $(X,T)$. We will make explicit  the direct link between the
labelling and the cocycle on the quotient.  In the next few sections, we
describe an ordered Bratteli diagram for the quotient of $(X_{\lambda},
T_{\lambda})$ by the $G$-action $\{\gamma_g\}_{g \in G}$.  The Bratteli
diagram thus constructed naturally maps to $(V,E, \geq)$ with unique
path lifting property.

Recall that by definition
$X_{\lambda} = \{x= (x_1, x_2, \cdots, x_n, \cdots ) \}$
where
\begin{enumerate}
\item[i)] $x_n = (x_{n,i})_{i \in \Z} \in ({\varpi}_{n, \lambda})^{\Z}$
\item[ii)] $j_{m,n} (x_{m,i}) =x_{n,i}$  for $m >n$ and
$i \in \Z$
\item[iii)] Given $n$ and $K$, $\exists \text{ } m >n$ and a vertex $v \in V_{m, \lambda}$
such that  the interval segment
$x_n [-K, K]:= (x_{n,-K}, x_{n,-K+1}, \cdots, x_{n,K})$ is obtained by
applying $j_{m,n}$ to an interval segment of the linearly ordered set of
paths from $v_0 $ to $v$.

\end{enumerate}

\noindent{\bf 2.11. Definition: ``Loops lift to loops''.}

If $(V, E, \geq)$ and $(V',E', \geq)$ are two ordered Bratteli diagrams
and $\pi: (V,E, \geq) \to (V',E', \geq)$ is a morphism with unique path-lifting
property we say that $\pi$ ``lifts loops to loops'' (or, that $\pi$ has the ``loops
lifting to loops" property) if the following
condition is satisfied:-

\noindent{\bf (``loops on the right") }Let $u' \in V'_m, v' \in V'_n$.  Let $k \leq m$ and $k \leq n$.  Suppose
$\alpha'$ and $\beta'$ are paths from $V'_k$  ranging at $u'$.  Suppose
$\gamma'$ and $\delta'$ are paths from $V'_k$ ranging at $v'$.  In the
lexicographic order induced on paths from $V'_k$ to $u'$ assume that
$\alpha'$ is the successor of $\beta'$.  Similarly for paths from $V'_k$
to $v'$ assume that $\gamma'$ is the successor of $\delta'$.  Assume
that $\alpha'$ and $\gamma'$  have the same source in $V'_k$  and that,
likewise, $\beta'$ and $\delta'$ have the same source in $V'_k$. (Thus, one has
a loop:- going from the range of  $\alpha'$ to the source of $\beta'$ via $\beta'$,
then to the range of $\delta'$ via  $\delta'$, then to the source of  $\alpha'$ via  $\gamma'$ and then to
the range of $\alpha'$ via  $\alpha'$.)

\noindent{\bf (`pull-back of the above loop to the left') }Now, let
$u \in V_m, v \in V_n$ and assume $\pi (u)=u', \pi (v)=v'$.  Let $\alpha
$ (resp. $\beta)$ be the unique path from $V_k$ to $V_m$ lying above
$\alpha'$, (resp. $\beta')$ and ranging at $u$.  Similarly, let $\gamma$
(resp. $\delta)$ be the unique path from $V_k$ to $V_n$ lying above
$\gamma'$ (resp. $\delta')$ and ranging at $v$.  If $\pi$ has the property
that under the above conditions,
$$\text{``\{source of }\beta=\text{ source of } \delta \}"
\Rightarrow ``\{\text{ source of }\alpha =\text{ source of }\gamma \}"$$
then we say that $\pi$ has the `loops lifting to loops' property.

%\newpage

\medskip
\centerline{\bf Diagram 3}
\begin{center}
\includegraphics[totalheight=8cm]{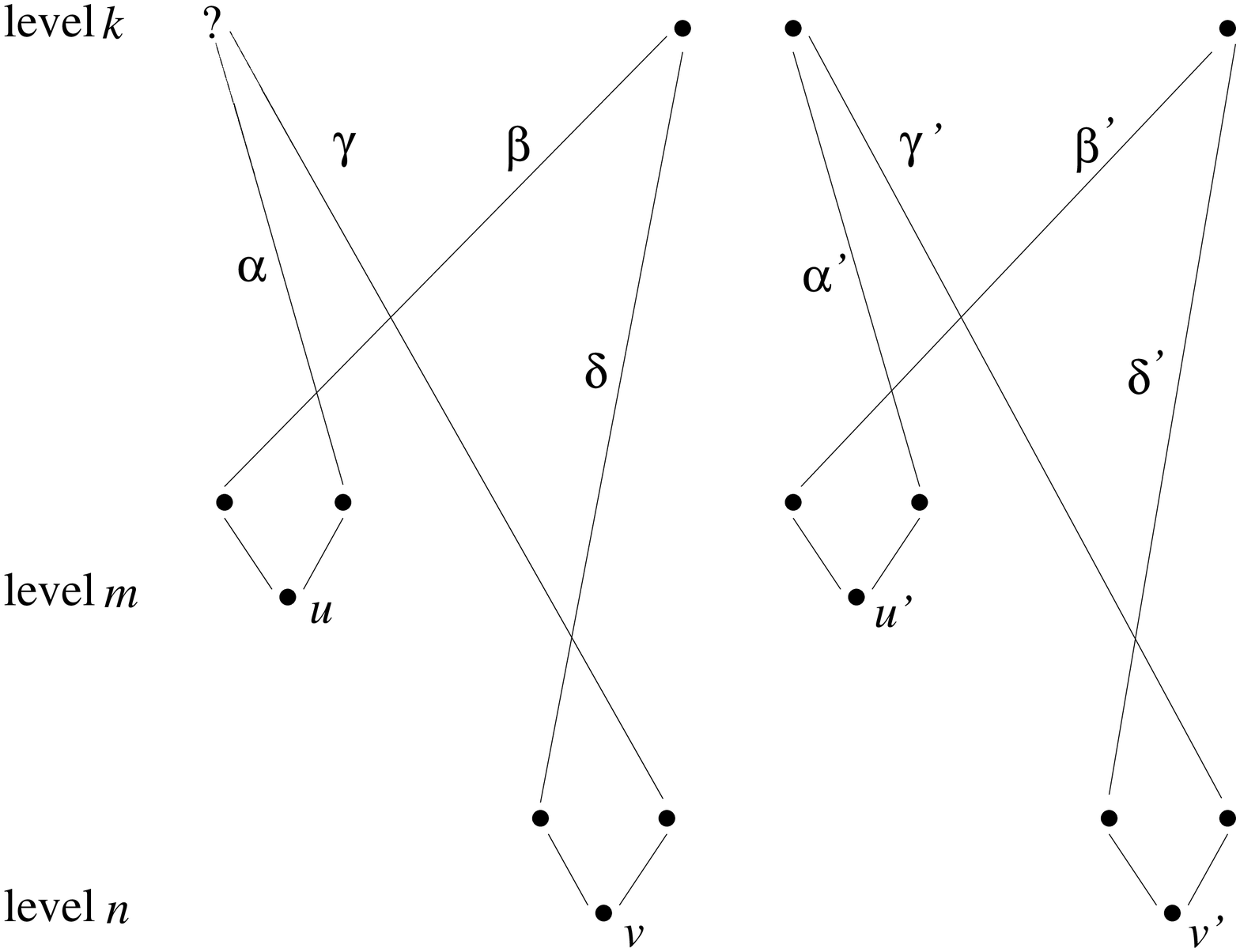} \end{center}

%\vskip7cm

Finally, we say that the labelling $\lambda:E\to G$ has the ``loops lifting to loops" property if
$(V_{\lambda},E_{\lambda},\geq) \to (V,E,\geq)$ has the above property.

\medskip
\noindent{\bf 2.11.1. Example. }Let $\beta :V\to G$ be a function. The labelling $\lambda$ defined on $E$ by
$\lambda (e) = \beta (r(e))^{-1}\beta (s(e))$ has the ``loops lifting to loops" property.\\
\noindent{\bf 2.11.2. Example. }Let $s,t\in G$. Let $\lambda$ be a $G$-valued stationary labelling on the
$2-adic$ odometer system({\it cf. }[GJ,pp.1691,1695]) whose image consists of $\{s,t\}$. 

The Bratteli diagram of the odometer has one vertex at each level and two edges, say, $min$ and $max$ from level $k$ to $k+1$
for $k \geq 1$, with $\lambda (min)=s$ and $\lambda (max)=t$.  In 2.11, we take \\
$\alpha' \text{ = } max$ from level 1 to 2,\\
$\beta'\text{ = } min$ from level 1 to 2,\\
$ \gamma'\text{ = } (min,max)$ from level 1 to 3,\\
$ \delta' \text{ = } (max,min)$ from level 1 to 3.\\
Denote by $w',u',v'$ the vertices at levels 1,2,3. Following the description of vertices and edges for $(V_{\lambda},
E_{\lambda})$ as in 2.5,\\ 
first (resp.second) edge ranging at $(u',g)$ has source $(w',gs)$ (resp.$(w',gt)$.\\ 
first (resp.second) edge ranging at $(v',h)$ has source $(u',hs)$ (resp.$(u',ht)$.\\
source $(w',gt)$ to range $(u',g)$ is a lift $\alpha$ of $\alpha'$.\\
source $(w',gs)$ to range $(u',g)$ is a lift $\beta $ of $\beta'$.\\
$(w',hts)$ to $(u',ht)$ to $(v',h)$ is a lift $\gamma$ of $\gamma'$.\\
$(w',hst)$ to $(u',hs)$ to $(v',h)$ is a lift $\delta$ of $\delta'$.\\
If the ``loops lifting to loops" property holds then for the above choice of $\alpha', \beta', \gamma', \delta'$
the condition that\\
$\text{    ``\{source of }\beta=\text{ source of } \delta \}"
\Rightarrow ``\{\text{ source of }\alpha =\text{ source of }\gamma \}"$ \\
yields the property in the group\\
$\text{    ``}\{(w',gs) = (w',hst)\}" \Rightarrow ``\{(w',gt) =(w',hts)\}".$\\
This reduces to the condition that $st^{-1}s^{-1}=ts^{-1}t^{-1}$.
One can easily see that analogous conditions with other choices of $\alpha', \beta', \gamma', \delta'$
are all consequences of the property $st^{-1}s^{-1}=ts^{-1}t^{-1}$.

Assume that $G$ is generated by $s,t$
and that $st^{-1}s^{-1}=ts^{-1}t^{-1}$. Then either $s=t$, (so $G$ is cyclic), or, $G$ has the following
simple description:

 Put $x=ts^{-1}$. Then, $x^2=ts^{-1}ts^{-1}=ts^{-1}ts^{-1}t^{-1}t=ts^{-1}st^{-1}s^{-1}t=s^{-1}t$.
 So, $t^{-1}xt=s^{-1}t=x^2$. Let $H$ be the subgroup of $G$ generated by $x$. Let $J$ be the
subgroup of $G$ generated by $t$. Note that elements of $J$ normalize $H$. Also, $J\cap H=\{e\}$ the
unit element. In fact since $x$ and $x^2$ are conjugate, the order of $x$ has to be $2k+1$ an odd
number. If $x^{\ell} = t^m \in J\cap H$, where, $0\leq \ell <2k+1$ then $x^{\ell}=t^m = t^{-1}x^{\ell}t=x^{2\ell}$,
which implies that $x^{\ell}=\{e\}$ which is not possible (unless $\ell =0$) since $\ell <$ order of $x$.  

 Hence, $G$ is the semi-direct
product of $J$ and $H$ for the action given by $t^{-1}xt=x^2$. The skew-product Bratteli diagram 
$(V_{\lambda},E_{\lambda},\geq)$ has the following connectivity structure.
Telescoping to levels
%$0=n_0,n_1,n_2,\cdots,n_i,\cdots$
$\{0,n_1=1,n_2,n_3,\cdots ,n_i,\cdots \}$ where $n_i=1+ (i-1) \times \{\text{order of } G\}$,
a vertex $(v,g)$ in level $n_i$ is connected
to a vertex $(w,g')$ in level $n_{i+1}$ by an edge if and only if $g \text{ and } g'$ belong to the same $H$-coset.

\smallskip
Our main result about the skew-product dynamical systems $(X_{\lambda},T_{\lambda})$ which were constructed in
$2.2$ are contained in the following.

\smallskip
\noindent{\bf 2.12.Theorem.} (I) {\it Suppose that the map $\pi :(V_\lambda,E_\lambda,\geq )\longrightarrow (V,E,\geq )$
has the property `loops lift to loops'. Then the quotient of $(X_\lambda , T_\lambda )$ by
the action of $G$ is canonically isomorphic to $(X,T)$.}

(II) {\it  In general, there is a commutative diagram}
$$\begin{matrix}
(V_\lambda, E_\lambda ,\geq )&\buildrel \pi \over \longrightarrow &(V,E,\geq )\\
\rho'\uparrow &&\uparrow \rho \\
(\widetilde V_\mu , \widetilde E_\mu ,\geq )&\buildrel {\widetilde \pi } \over \longrightarrow
&(\widetilde V,\widetilde E,\geq )
\end{matrix}$$
{\it and the induced diagram}
$$\begin{matrix}
(X_\lambda, T_\lambda   )&\buildrel {\pi_B} \over \longrightarrow &(X,T  )\\
\rho'_B\uparrow &&\uparrow \rho_B \\
(\widetilde X_\mu , \widetilde T_\mu   )&\buildrel {{\widetilde \pi }_B} \over \longrightarrow
&(\widetilde X,\widetilde T  )
\end{matrix}$$
{\it between the corresponding dynamical systems, where}

 (i) {\it  All the horizontal and vertical arrows in the first diagram have the `unique path lifting property'},

(ii) {\it $(\widetilde V_\mu , \widetilde E_\mu ,\geq )$ is the skew-product Bratteli diagram associated
to the $G$-valued labelling $\mu $ on  $ (\widetilde V,\widetilde E,\geq )$ defined by   $\mu = \lambda \circ \rho$}

(iii) {\it $\rho'$ commutes with the $G$-action and $\rho_B' $ is an isomorphism and}

(iv) {\it the map $\widetilde \pi  : (\widetilde V_\mu , \widetilde E_\mu ,\geq )   \longrightarrow
(\widetilde V,\widetilde E,\geq )$ has the property ``loops lift to loops".} ({\it Consequently, by Part} (I),
{\it the quotient of $(X_\lambda , T_\lambda )$ by
the action of $G$ is isomorphic to $(\widetilde X,\widetilde T)$}.)

\noindent{\bf Remark. }Note that Part (II) is almost a converse to Part (I).
If we are given that
the quotient of $(X_\lambda , T_\lambda )$ by
the action of $G$ is isomorphic to $(X,T)$, we are not deducing from this that
$\pi :(V_\lambda,E_\lambda,\geq )\longrightarrow (V,E,\geq )$
has the property `loops lift to loops'.  Instead, Part (II) gives us a
closely related morphism with the property `loops lift to loops' which can be
nicely fit into a square with the given morphism $\pi$ as one of its sides and
in which all the arrows have the unique path lifting property; the
vertical arrows are compatible with the given labelling and its pull-back and
induce isomorphisms in the corresponding dynamical systems.

\noindent{\bf Poof of part (I).} By hypothesis $\pi :(V_\lambda,E_\lambda,\geq )\longrightarrow  (V,E,\geq )$ has the property `loops lifts
 to loops'. Suppose that $\pi_B :(X_\lambda,E_\lambda,\geq )\longrightarrow  (X,E,\geq )$ is not  an isomorphism, modulo the
$G$-action on $X_\lambda $. Let $y$ and $z$ be two different $G$-orbits in $X_\lambda $, which have the same image $x$
in $X$. We employ the notation  of $(2.9)$  and $(2.10)$. Write

$y=(y_n)_n,\hskip0.5cm y_n=(y_{n,i})_{i\in \Z},\hskip0.5cm y_{n,i}\in \varpi_{n,\lambda } $

$z=(z_n)_n,\hskip0.5cm z_n=(z_{n,i})_{i\in \Z},\hskip0.5cm z_{n,i}\in \varpi_{n,\lambda } $

$x=(x_n)_n,\hskip0.5cm x_n=(x_{n,i})_{i\in \Z},\hskip0.5cm x_{n,i}\in \varpi_n $.

\noindent By the hypothesis, $y$ and $z$ have the image $x$ in $X$ and also the $G$-orbits of $y$ and
$z$ are different. By shifting $y$ and $z$ by a suitable power of the shift operator $T_\lambda $ and replacing
$z$ if necessary by another element within the orbit $\{ \gamma_g\cdot z\mid g\in G\} $, we can assume that, for some
$k\in \N$, $y_{k,0}=z_{k,0}$ and $y_{k,1}\neq z_{k,1}$. In view of (1.8) it follows that\\
{\bf (2.12.1.)}  $y_{k,0}$ is the top floor of a
$\varpi_{k,\lambda }$-tower $\varpi (u,g_1)$ \\
and \\
{\bf (2.12.2.)}  $y_{k,1}$ is the lowest floor of a
$\varpi_{k,\lambda }$-tower $\varpi (w,g_2)$.\\
In the same way, \\
{\bf (2.12.3.)}  $z_{k,0}$ is the top floor of the
$\varpi_{k,\lambda}$-tower $\varpi (u,g_1)$\\
and \\
{\bf (2.12.4.)}  $z_{k,1}$ the lowest floor of a
$\varpi_{k,\lambda }$-tower $\varpi (w ,g_3)$.\\
By property ((1.10), (iii)) there exists $m\geq k$, a vertex $(a,h)\in V_{m,\lambda}$ and two successive floors
$F^0,F^1$ of the $\varpi_{m,\lambda}$ tower $\varpi (a,h)$ such that $j_{m,k}(F^0)=y_{k,0}$ and $j_{m,k}(F^1)=y_{k,1}$.
Similarly, there exists $n\geq k$, a vertex
$(b,h^{\tilde{}})\in V_{n,\lambda}$
and two successive floors
$\widetilde F^0,\widetilde F^1$ of the $\varpi_{n,\lambda}$-tower $\varpi (b,h^{\tilde{}})$ such that $j_{m,k}(\widetilde F^0)=y_{k,0}$ and $j_{m,k}(\widetilde F^1)=y_{k,1}$.
The floors $F^0$ and $F^1$ correspond to two successive paths from $\{ \ast \} \in V_{0,\lambda}$ to
$(a,h)\in V_{m,\lambda}$ in the linearly ordered set of paths from $\{ \ast \} $ to $(a,h)$.
The path $F^0$ is traced by first tracing a path $F^0[0,k]$ from level $0$ to level $k$ and following it by a path
from level $k$ to level $m$. The part $F^0[0,k]$ being the truncation $j_{m,k}(F^0)$ represents the floor $y_{k,0}$.
By $(2.12.1)$, $y_{k,0}$ is the unique maximal path from $\{ \ast \} \in V_{0,\lambda}$ to
$(u,g_1)\in V_{k,\lambda}$. Similarly, the path $F^1$ from level $0$ to level $m$ initially traces  the truncation
$F^1[0,k]=j_{m,k}(F^1)$ which represents the floor $y_{k,1}$  followed by $F^1[k,m]$ from
level $k$ to level $m$. By $(2.12.2)$, $y_{k,1}$ is the unique minimal path (lowest floor) from
$\{ \ast \}\in V_{0,\lambda}$ to $(w,g_2)\in V_{k,\lambda}$. Since $F^1$ is the successor of $F^0$ in the
lexicographic order of paths from $\{ \ast \} $ to $(a,h)\in V_{m,\lambda}$, we conclude from the
above observation that

{\bf (2.12.5)} $F^1[k,m]$ is the successor of $F^0[k,m]$ in the lexicographic order of paths from $V_{k,\lambda}$
to $(a,h)\in V_{m,\lambda}$.

Similarly, if $\widetilde F^0[k,n]$  and $\widetilde F^1[k,n]$ represent the truncation
of $\widetilde F^0$ and $\widetilde F^1$ respectively from level $k$ to level $n$ then

{\bf (2.12.6)}  $\widetilde F^1[k,n]$ is the successor of $\widetilde F^0[k,n]$ in the lexicographic order of paths from $V_{k,\lambda}$
to $(a,h^{\tilde{}})\in V_{n,\lambda}$.

We show how this leads to a contradiction of the property `loops lift to loops'. Denote the paths
$F^0[k,m]$, $F^1[k,m]$, $\widetilde F^0[k,n]$ and $\widetilde F^1[k,n]$ by $\beta , \alpha , \delta $ and $\gamma $
respectively. Let their images under $\pi $ be denoted by $\beta' , \alpha' , \delta' $ and $\gamma' $ respectively.

We will observe that with the above data we have something on the left which is almost a loop, but not actually a loop
eventhough its image on the right is a loop ({\it cf. } (2.11)).

\noindent{\bf (object on the right which is a loop) }
$\alpha'$ and $\beta'$ are paths from $V_k$  ranging at $a$.  Also
$\gamma'$ and $\delta'$ are paths from $V_k$ ranging at $b$.
$\alpha'$ is the successor of $\beta'$.  For paths from $V_k$
to $V_n, \gamma'$ is the successor of $\delta'$. The paths $\alpha'$ and $\gamma'$  have the
same source in $V_k$, $\beta'$ and $\delta'$ have the same source in $V_k$. (Thus, one has
a loop:- going from the range of  $\alpha'$ to the source of $\beta'$ via $\beta'$,
then to the range of $\delta'$ via  $\delta'$, then to the source of  $\alpha'$ via  $\gamma'$ and then to
the range of $\alpha'$ via  $\alpha'$.)

\noindent{\bf (object on the left lying over the above loop)} $(a,h )\in V_{m,\lambda}, (b,h^{\tilde{}})\in V_{n,\lambda}$,  $\pi (a,h )=a, \pi (b,h^{\tilde{}})=b$.
 $\alpha
$ (resp. $\beta)$ is  the unique path from $V_{k,\lambda}$ to $V_{m,\lambda}$ lying above
$\alpha'$, (resp. $\beta')$ and ranging at $(a,h)$.  Similarly,  $\gamma$
(resp. $\delta)$ is the unique path from $V_{k,\lambda}$ to $V_{n,\lambda}$, lying above
$\gamma'$ (resp. $\delta')$ and ranging at $(b,h^{\tilde{ }})$. $\alpha$ is the successor of $\beta$.
$\gamma$ is the successor of $\delta$.

\noindent{\bf (object on the left is almost a loop) } Moreover  source of $\beta=$ source of $\delta $.

But  source of $\alpha \neq $ source of $\gamma $,
which contradicts the `loops lifting to loops' property. \hfill $\square $

Now we begin the proof of the part II in the statement of the theorem.

\noindent{\bf 2.13.} We will now define two nested sequences of K-R partitions of
$X_{\lambda}$, both acted on by the group $G$.  Recall $V_{n, \lambda} =
V_n \times G$ (for $n \geq 1)$.  For $(v,g) \in V_n \times G$,
let $y$ be a path from $\{ *\}$ to $(v,g)$ in $(V_{\lambda}, E_{\lambda},
\geq)$.  So, $y$ is a  `floor'  belonging to the $\varpi_{n, \lambda}$-
tower $\varpi (v,g)$ parametrized by $(v,g)$.
Put ${\mathcal F}_y = \{x= (x_1, x_2, \cdots, x_n, \cdots ) \in X_{\lambda}  \mid
x_{n,0} =y \}$
$${\mathcal P}_{n, \lambda} = \{{\mathcal F}_y \mid y \in \varpi (v,g), (v,g) \in V_{n, \lambda}\}.$$
Then $\{{\mathcal P}_{n, \lambda} \}_n$ is a nested sequence of K-R partitions of $X_{\lambda}$
acted on by $ G$.  But, the topology of $X_{\lambda}$ need not be spanned by
the collection of clopen sets $\{{\mathcal F}_y\}, (y \in \varpi (v,g), (v,g)
\in V_{n, \lambda}, n \in \N)$.  In contrast, the topology of $X_{\lambda}$
is indeed spanned by the collection of clopen sets in another nested sequence
$\{{\mathcal Q}_{n, \lambda}\}_n$ of K-R partitions, defined below.
Let $\varpi = \varpi (u,g_1), \varpi' = \varpi (v,g_2),
\varpi'' = \varpi (w,g_3)$
be three $\varpi_{n, \lambda}$-towers and $y$ a floor of $\varpi'$.
For any $x \in X_{\lambda}$ and for any $n$ if $x_{n,i}$ is a floor of a
$\varpi_{n, \lambda}$-tower $\overline{\varpi}$, then for some
$a,b\in \Z$ such that $a \leq i \leq b$, the segment $x_n [a,b]$ is
just the sequence of floors in $\overline{\varpi}$.  We define

${\mathcal F}(\varpi, \varpi', \varpi''; y) =$ the clopen subset of
${\mathcal F}_y$ consisting of the elements $x= (x_1, x_2, \cdots, x_n, \cdots)$ with the property that for some $a_1 < a_2 \leq 0 < a_3 < a_4 \in \Z$, the segment $x_n [a_1,
a_2-1]$ is the sequence of floors of $\varpi$, the segment $x_n
[a_2, a_3 -1]$ is  the sequence of floors of $\varpi'$ and the segment
$x_n [a_3, a_4]$ is the sequence of floors of $\varpi''$.
Some of the sets  ${\mathcal F}(\varpi, \varpi', \varpi'';y)$ may be empty, but
the non-empty
sets ${\mathcal F}(\varpi, \varpi', \varpi'';y)$ form a K-R partition which we denote by
${\mathcal Q}_{n,\lambda}$. For fixed $\varpi, \varpi', \varpi'' $
the subcollection $\{{\mathcal F}(\varpi,
\varpi', \varpi'';y)\}$ as $y$ varies through  the floors of
$\varpi'$, is a ${\mathcal Q}_{n, \lambda}$-tower parametrized by
$[(u,g_1), (v,g_2),(w, g_3)]$.  We denote this ${\mathcal Q}_{n, \lambda}$-tower by
${\mathcal S}_{(\varpi,
\varpi', \varpi'')}$.  The floors of the tower ${\mathcal S}_{(\varpi, \varpi', \varpi'')}$
are $\{{\mathcal F}(\varpi, \varpi', \varpi'';y)\}$ as $y$
runs through the sequence of floors of $\varpi'$.

\noindent{\bf 2.14. Lemma.}  {\it $\{{\mathcal Q}_{n, \lambda}\}_n$ is a nested sequence of
K-R partitions of $X_{\lambda}$ acted on by $G$.  The topology of $X_{\lambda}$ is spanned
by the clopen sets in this sequence of partitions.}

\paragraph*{Proof.} Each ${\mathcal Q}_{n, \lambda} = \{{\mathcal F}(\varpi, \varpi', \varpi'';y)\} $ (for a fixed $n)$ is a K-R
partition. This is evident on going through the definitions.  Since the
infimum of the height of ${\varpi}_{n, \lambda}$-towers $(\varpi, \varpi'$ etc.)
goes to infinity as $n \to \infty$, it is evident from the definition of the
topology of $X_{\lambda}$ that the clopen sets ${\mathcal F} (\varpi, \varpi',
\varpi'';y)$, (for various $n, {\varpi}_{n,\lambda}$- towers  $\varpi,\varpi',\varpi"$, floors
$y$ of $\varpi')$ span the topology of $X_{\lambda}$. \hfill $\square $

\noindent{\bf 2.15. Remark.} The above construction involving `triples' can be
introduced starting with any simple ordered Bratteli diagram, (not necessarily
the diagram $(V_{\lambda}, E_{\lambda}, \geq)$ which had a  $G$-action).  For this, let
$(V, E, \geq)$ be an arbitrary simple, ordered Bratteli diagram.  Define
$(V^{\mathcal Q}, E^{\mathcal Q}, \geq )$ as follows:  $V^{\mathcal Q}_0 =\{*\}$, a
single point.

$V^{\mathcal Q}_n$ consists of triples $(u,v,w) \in V_n\times V_n \times V_n$ such that for
some $y \in V_m$ where $m>n$, the level-$m$ tower ${\varpi}_y$ passes successively
through the level-$n$ tower ${\varpi}_u$, then ${\varpi}_v$ and then ${\varpi}_w$.
An edge $ \tilde{e} \in  E^{\mathcal Q}_n$ is a triple $(u,e,w)$ such that
$e$ is an edge of $(V,E) $ and $(u,r (e),w) \in V^{\mathcal Q}_n$.  Let
\begin{enumerate}
\item[]$\{e_1, e_2, \cdots, e_k\}$ be all  the edges  in  $r^{-1} (r(e))$,
\item[]$\{f_1, f_2, \cdots, f_{\ell}\}$ be all the edges in  $r^{-1} (u)$ and
\item[]$\{g_1, g_2, \cdots, g_m \}$ be all the edges in  $r^{-1} (w)$.
\end{enumerate}
The sources of $(u,e_1, w), (u, e_2, w), \cdots, (u,e_k, w)$ are defined
to be \\
$ (s(f_{\ell}), s(e_1),s(e_2)), (s(e_1),s (e_2), s(e_3)),
\cdots, (s(e_{k-1}), s(e_k), s(g_1))$ \\
respectively. The range of
$(u,e,w)$ is of course $(u,r (e), w)$.

The map $(u,v,w) \mapsto v, (u,e,w) \mapsto e$ from $(V^{\mathcal Q}, E^{\mathcal Q})$
to $(V,E)$ has unique path lifting property; in particular it gives rise to
the ordered Bratteli diagram $(V^{\mathcal Q}, E^{\mathcal Q}, \geq )$.

\noindent{\bf (Example (2.6) continued). }We illustrate the tripling construction
described above for the substitutional system $A=\{X,Y\},\sigma (X)=XXY,
\sigma (Y)=XYY$.
One can easily verify that the tripling yields a
stationary ordered Bratteli diagram corresponding to substitutional system with
alphabet set $A^{\mathcal Q} \subseteq A \times A \times A$ described by
$$A^{\mathcal Q} =\{[X,X,Y], [Y,X,X], [X,Y,X], [X,Y,Y], [Y,X,Y], [Y,Y,X]\}$$
and substitution $\sigma^{\mathcal Q}$ described by
$$
\begin{array}{lll}
\sigma^{\mathcal Q} ([X,X,Y]) &=& [Y,X,X][X,X,Y][X,Y,X] \\
\sigma^{\mathcal Q} ([Y,X,X]) &=& [Y,X,X][X,X,Y][X,Y,X] \\
\sigma^{\mathcal Q} ([X,Y,X]) &=& [Y,X,Y][X,Y,Y][Y,Y,X] \\
\sigma^{\mathcal Q} ([X,Y,Y]) &=& [Y,X,Y][X,Y,Y][Y,Y,X] \\
\sigma^{\mathcal Q} ([Y,X,Y]) &=& [Y,X,X][X,X,Y][X,Y,X] \\
\sigma^{\mathcal Q} ([Y,Y,X]) &=& [Y,X,Y][X,Y,Y][Y,Y,X]
\end{array}
$$

In this example, since $\sigma$ is a proper  substitution, $\sigma^{\mathcal Q}$
depends only on the middle coordinate.  Later we will illustrate with a more
interesting example.

\noindent{\bf 2.16. The ordered Bratteli diagrams $B^{\mathcal Q}(\lambda)$
and  $\overline{B}^{\mathcal Q}(\lambda)$. }
The construction $(1.6)$ of an
ordered Bratteli diagram applied to the  nested sequence of K-R partitions
 $\{{\mathcal Q}_{n, \lambda}\}_n$  gives rise
to an ordered Bratteli diagram.  The vertices $V_{n, \lambda}^{{\mathcal Q}}$ are in 1-1
correspondence with the towers ${\mathcal S}_{[\varpi (u,g_1), \varpi (v,g_2),
\varpi (w,g_3)]}$ of ${\mathcal Q}_{n, \lambda}$.  Note that not all choices 
$[(u,g_1), (v, g_2),\\ (w,g_3)]$ may give rise to a non-empty set ${\mathcal S}_{[\varpi (u,g_1),
\varpi (v,g_2), \varpi (w, g_3)]}$.  We may denote the vertex corresponding to (non-empty)
${\mathcal S}_{[\varpi (u,g_1), \varpi  (v,g_2), \varpi (w,g_3)]}$ by  $[(u,g_1), (v,g_2), (w,g_3)]$.  
If $[(u,g_1), (v,g_2), (w,g_3)] \in
V_{n, \lambda}^{\mathcal Q}$, then $[(u,gg_1),\\ (v,gg_2), (w,gg_3)] \in V_{n, \lambda}^{\mathcal Q}, \forall g \in G$. 
So $G$ acts on $V_{n, \lambda}^{Q}$ freely if
$n \geq 1$.  As in $(1.6)$ the edges connecting $V_{n, \lambda}^{\mathcal Q}$
and $ V_{n+1, \lambda}^{\mathcal Q}$ and the linear order between edges with the same range 
simply reflects the sequence of ${\mathcal Q}_{n, \lambda}$-towers traversed by a
given ${\mathcal Q}_{n+1, \lambda}$- tower.  It is then clear that $G$-acts on the
ordered Bratteli diagram $(V_{\lambda}^{\mathcal Q}, E_{\lambda}^{\mathcal Q}, \geq )$.  We
denote by $(\overline{V}_{\lambda}^{\mathcal Q}, \overline{E}_{\lambda}^{\mathcal Q}, \geq )$
the ordered Bratteli diagram which is the quotient by $G$-action on
$(V_{\lambda}^{\mathcal Q}, E_{\lambda}^{\mathcal Q}, \geq)$. We write $B^{\mathcal Q}(\lambda)$
and  $\overline{B}^{\mathcal Q}(\lambda)$
for the Bratteli diagrams $(V_{\lambda}^{\mathcal Q}, E_{\lambda}^{\mathcal Q}, \geq)$ and
$(\overline{V}_{\lambda}^{\mathcal Q}, \overline{E}_{\lambda}^{\mathcal Q}, \geq )$ respectively.

\noindent{\bf 2.17. Examples 2.4 and 2.7 continued.} For the ordered Bratteli diagram $(V_{\lambda},E_{\lambda},\geq)$
which arose in (2.7) from the labelling $\lambda$ defined there,
the `tripling' construction $(2.15)$
yields a Bratteli diagram $(V^{\mathcal Q}_{\lambda}, E^{\mathcal Q}_{\lambda}, \geq )$,  which is
associated to the primitive substitutional system
below.  The alphabet set $A_{\lambda}^{\mathcal Q}$ consists of the twenty alphabets
\begin{enumerate}
\item[] [x,x,y], [x,y,x], [x,y,x'], [y,x,x], [y',x,x]
\item[] [y,x,y'], [y',x,y'],[y',y,x'], [y',y,x], [x',y,y']
\item[] [x',x',y'], [x',y',x'], [x',y',x], [y',x',x'], [y,x',x']
\item[] [y',x',y], [y,x',y],[y,y',x], [y,y',x'], [x,y',y]
\end{enumerate}
The substitution $\sigma_\lambda^{\mathcal Q}:A_\lambda^{\mathcal Q}\longrightarrow (A_\lambda^{\mathcal Q})^+$ is given by

$\sigma_\lambda^{\mathcal Q} ([x,x,y])=[y,x,x][x,x,y][x,y,x]$

$\sigma_\lambda^{\mathcal Q} ([x,y,x])=[y,x,y'][x,y',y][y',y,x]$

$\sigma_\lambda^{\mathcal Q} ([x,y,x'])=[y,x,y'][x,y',y][y',y,x]$

$\sigma_\lambda^{\mathcal Q} ([y,x,x])=[y,x,x][x,x,y][x,y,x]$

$\sigma_\lambda^{\mathcal Q} ([y',x,x])=[y',x,x][x,x,y][x,y,x]$

$\cdots $ {\it etc..}

\noindent The general rule for the substitution $\sigma_\lambda^{\mathcal Q}$ is expalined in the more general
 situation in $(2.18)$, below. The group $\Z_2$ acts on $(A_{\lambda}^{\mathcal Q}, \sigma_{\lambda}^{\mathcal Q})$ and
on $(V_{\lambda}^{\mathcal Q},E_{\lambda}^{\mathcal Q},\geq)$ via the transposition $x\rightleftarrows x', y\rightleftarrows y'$.
The quotient substitutional system $(\overline{A}_{\lambda}^{\mathcal Q}, \overline{\sigma}_{\lambda}^{\mathcal Q})$ with
ten alphabets maps to $(A,\sigma)$ of example (2.6). The induced map between the dynamical systems is not
an isomorphism.

\noindent{\bf 2.18. Tripling for a substitutional system $(A, \sigma)$.} If $(A, \sigma)$ is a substitutional system, define $A^{\mathcal Q}$ to be the
unique smallest subset of $A \times A \times A$ with the property that
$(a,b,c) \in A^{\mathcal Q}$ if and only if the word $abc$ occurs as a subword of
$\sigma^n (d)$ for some $d \in A$ and some $n$.  Define
$$\sigma^{\mathcal Q}: A^{\mathcal Q} \to (A^{\mathcal Q})^+$$
by $\sigma^{\mathcal Q} [(a,b,c)] =(a_m, b_1, b_2) \cdot (b_1,b_2,b_3) \cdots
(b_{n-2}, b_{n-1}, b_n) \cdot (b_{n-1}, b_n, c_1)$,  where $\sigma (b) =b_1 \cdot b_2 \cdots b_n$, and $a_m$ is the last alphabet in $\sigma (a)$, while $c_1$ is the first alphabet in $\sigma (c)$.
Suppose moreover that the ordered Bratteli diagram associated to $(A,\sigma)$ is equipped with a stationary labelling $\lambda$.
Then following the construction of $(2.5)$ we get a skew-product substitutional system $(A_{\lambda},\sigma_{\lambda})$.
The tripling described above applied to this $(A_{\lambda},\sigma_{\lambda})$ gives a substitutional system $(A^{\mathcal Q}_\lambda,\sigma^{\mathcal Q}_\lambda)$.

\noindent{\bf 2.19. }Even if we start with a primitive ({\it cf. }[DHS,3.3.1]) and proper substitutional system and a stationary labelling on the edges with
the additional property that the maximal and minimal edges are both labelled by the trivial element of the group,
the property ``loops lift to loops" may generally fail; the quotient of the dynamical system associated to
$(A_{\lambda}^{\mathcal Q}, \sigma_{\lambda}^{\mathcal Q})$
by the $G$-action will in general be a nontrivial extension of the dynamical system of $(A,\sigma)$.
This is an important difference between our construction of
skew-product systems arising from $G$-valued labellings and
Matui's construction of skew-product systems arising from $G$-valued cocycles.

\noindent{\bf 2.20.} It now remains to finish the proof of part (II) of Theorem (2.12). Our construction of
$\overline{B}^{\mathcal Q}(\lambda)=\left(  \overline{V}_\lambda^{\mathcal Q},\overline{E}_\lambda^{\mathcal Q},\geq \right) $
in $(2.16)$
was motivated precisely to serve as a candidate for the $\widetilde B$ in the statement of part II
of the theorem. Thus, we set $\widetilde{V}=\overline{V}_\lambda^{\mathcal Q}$,  $\widetilde{E}=\overline{E}_\lambda^{\mathcal Q}$,
$\widetilde{B}=(\widetilde{V},\widetilde{E},\geq )=\left( \overline{V}_\lambda^{\mathcal Q},\overline{E}_\lambda^{\mathcal Q},\geq \right)
=\overline{B}^{\mathcal Q}(\lambda)$. A vertex of $\overline{V}_{n,\lambda }^{\mathcal Q}$ is represented by the
$G$-orbit of a triple $[(u,g_1),(v,g_2),(w,g_3)]$. Define $\rho :\overline{V}_\lambda^{\mathcal Q}\longrightarrow V$
by sending the above vertex to $v\in V_n$. The set $\overline{E}_{n,\lambda }^{\mathcal Q}$
of edges of $\overline{B}^{\mathcal Q}(\lambda)$ from $\overline{V}_{n-1,\lambda }^{\mathcal Q}$
to $\overline{V}_{n,\lambda }^{\mathcal Q}$ ranging at $[(u,g_1),(v,g_2),(w,g_3)]$ is represented by the triple
$[(u,g_1),(e,g_2),(w,g_3)]$ where $e$ is an edge from $V_{n-1}$ to $V_n$ ranging at $v$. Define
$\rho : \overline{E}_\lambda^{\mathcal Q}\longrightarrow E$ by sending the above edge (namely the $G$-orbit of
$[(u,g_1),(e,g_2),(w,g_3)]$) to $e$.

Define  a labelling $\mu $ on $\overline{E}_\lambda^{\mathcal Q}$ by $\mu = \lambda \circ \rho $. We show that the
corresponding skew-product of $\left( \overline{V}_\lambda^{\mathcal Q},\overline{E}_\lambda^{\mathcal Q},\geq \right) $
by $\mu $ is precisely $(V_\lambda^{\mathcal Q},E_\lambda^{\mathcal Q},\geq )$. For this we define a map
$$\Phi :(V_\lambda^{\mathcal Q},E_\lambda^{\mathcal Q},\geq )\longrightarrow
(\overline{V}_\lambda^{\mathcal Q}\times G,\overline{E}_\lambda^{\mathcal Q}\times G,\geq )$$
as follows: let $x=[(u,g_1),(v,g_2),(w,g_3)]\in V_\lambda^{\mathcal Q}$. Let the $G$-orbit of
$x$ be denoted by $\overline x\in \overline{V}_\lambda^{\mathcal Q}$. Then, define
$\Phi (x)=(\overline x,g_2)$. Similarly, let $\eta =[(u,g_1),(e,g_2),(w,g_3)]$ be an edge of $E_\lambda^{\mathcal Q}$.
Let the $G$-orbit of $\eta $ be denoted by $\overline \eta \in \overline{E}_\lambda^{\mathcal Q}$. Then, define
$\Phi (\eta )=(\overline \eta ,g_2)$. Then $\Phi $ is an isomorphism between the two ordered Bratteli diagrams.

Referring to the requirements in the statement of Part (II) of Theorem $(2.12)$,
the map $\rho':(V_\lambda^{\mathcal Q},E_\lambda^{\mathcal Q},\geq )\longrightarrow (V_\lambda , E_\lambda ,\geq )$
is easy to define:
$$\rho'[(u,g_1),(v,g_2),(w,g_3)]=(v,g_2)\hskip0.5cm \hbox{and}\hskip0.5cm
\rho'[(u,g_1),(e,g_2),(w,g_3)]=(e,g_2).$$
The map $\rho':B^{\mathcal Q}(\lambda) = (V_\lambda^{\mathcal Q},E_\lambda^{\mathcal Q},\geq )
\longrightarrow (V_\lambda , E_\lambda ,\geq )=B(\lambda ) $ induces a map $\rho'_B$ between the corresponding
dynamical systems $X_{B^{\mathcal Q}(\lambda )}$ and $X_{B(\lambda )}$. The inverse $\beta $ of this map is defined as follows:

Let $z\in X_{B(\lambda)} $. So $z=(z_1,\cdots , z_n,\cdots )$, where $z_n=(z_{n,i})_{i\in \Z}\in (\varpi_{n,\lambda })^{\Z}$.
Fix $n$ and $i$. So $z_{n,i}$ is a floor of  a $\varpi_{n,\lambda }$-tower ${\mathcal S}_{(v,g_2)}$. Choose $a,b,c,d\in \Z$ such that
$a<b\leq i<c<d$, the interval segment $z_n[b,c-1]$ is the sequence of all the floors of the tower ${\mathcal S}_{(v,g_2)}$,
$z_n[a,b-1]$ is the sequence of all the floors of some $\varpi_{n,\lambda }$-tower ${\mathcal S}_{(u,g_1)}$ and $z_n[c,d-1]$ is the sequnce of all the floors of some $\varpi_{n,\lambda }$-tower ${\mathcal S}_{(w,g_3)}$. Define
$$y=(y_1,\cdots ,y_n,\cdots ) \in X_{B^{\mathcal Q}(\lambda )} \hskip0.7cm \hbox{where } y_n=(y_{n,i})_{i\in \Z}\in \left( \varpi_{n,\lambda }^{\mathcal Q}\right)^{\Z}$$
by $y_{n,i}=((u,g_1),(v,g_2),(w,g_3);z_{n,i})$. This map defines an inverse of $\rho'_B$.

Finally, we prove the property `loops lift to loops' for the map
$$(V_\lambda^{\mathcal Q}, E_\lambda^{\mathcal Q},\geq )\longrightarrow
(\overline{V}_\lambda^{\mathcal Q},\overline{E}_\lambda^{\mathcal Q},\geq ).$$
To verify the property `loops lift to loops' for the above map, we start with some data on the left-side
which is `almost' a loop, we make a further assumption that the image is actually a loop and then
we have to conclude that the data we started with on the left-side must in fact be a loop.

\noindent{\bf 2.21. The data on the left which is {\it almost} a loop.} Let $\hat{u}\in V_{m,\lambda}^{\mathcal Q}$ and $\hat{v}\in V_{n,\lambda }^{\mathcal Q}$.
Let $k\leq m$ and $k\leq n$.
Let $\alpha $  and $\beta $ be
paths from $V_{k,\lambda }^{\mathcal Q}$ ranging at $\hat{u}$. Let $\gamma ,\delta $ be paths from $V_{k,\lambda }^{\mathcal Q}$ ranging at $\hat{v}$.
Assume that $\beta$ and $\delta$ have the same source in
$ V_{k,\lambda }^{\mathcal Q}$.

\noindent{\bf 2.22. The data `image on the right is {\it actually} a loop'.} Denote by $\alpha', \beta',\gamma',\delta'$ the
paths in $(\overline{V}_\lambda^{\mathcal Q},\overline{E}_\lambda^{\mathcal Q},\geq )$ which are the
images of $\alpha,\beta,\gamma,\delta$.
Suppose that $\alpha'$ is the successor of $\beta'$ and that $\gamma'$ is the successor of $\delta'$. Assume that
$\alpha'$ and $\gamma'$
have the same source in $ \overline{V}_{k,\lambda}^{\mathcal Q}$    and that $\beta'$ and $\delta'$ have the same source in
$ \overline{V}_{k,\lambda}^{\mathcal Q}$.

To conclude that the data on the left must in fact be a loop we have to show that
$\alpha$ and $\gamma$
have the same source in $ V_{k,\lambda}^{\mathcal Q}$.

Let source of $\beta =$ source of $\delta =  [(a,g_1),(b,g_2),(c,g_3)]\in V_{m,\lambda}$. Since $\alpha$ is the succssor of $\beta$,
necessarily the source of $\alpha$ has the form  $[(b,g_2),(c,g_3),(d,g_4)]$ for some $d \in V_k$ and $g_4 \in G$.
Likewise, since $\gamma$ is the successor of $\delta$ the source of $\gamma$ must be of the form
$[(b,g_2),(c,g_3),(d',g_5)]$, for some $d' \in V_k$ and $g_5 \in G$. But it is a consequence of the assumption ``source of $\alpha'$
= source of $\gamma'$" that $[(b,g_2),(c,g_3),(d,g_4)]$ and $[(b,g_2),(c,g_3),(d',g_5)]$ must lie in the same $G$-orbit,
i.e., for some $h\in G$,
$[(b,g_2),(c,g_3),(d,g_4)]$ = $[(b,hg_2),(c,hg_3),(d',hg_5)]$.

\noindent{}Therefore, we conclude that $d=d'$ and $g_4= g_5$.  Thus, the property ``loops lift to loops" is
verified. 

 This ends the proof of the theorem.\hfill $\square $

\smallskip

\noindent{\bf 2.23. }We already gave an example ({\it cf. }2.17) where the map $\rho_B$ of Theorem 2.12 is not an
isomorphism. To give another example let $G=\mathfrak{S}_3$  be the group consisting of the six permutations of the
symbols $\{1,2,3\}$. Let $s$ and $t$ be the elements of $\mathfrak{S}_3$ corresponding to the 2-cycle $(1,2)$ and
the 3-cycle $(1,2,3)$ respectively. Let $A$ be any alphabet set $\{a_{\tau} \mid \tau \in \mathfrak{S}_3 \}$.
Define a substitution $\sigma : A \to A$ by  $\sigma (a_{\tau})=a_{\tau t} a_{\tau s}$. The group $\mathfrak{S}_3$ acts
on $(A,\sigma)$ by $\gamma_g(a_{\tau})=a_{g\tau}$. This is a skew-product system  ({\it cf. }2.5) associated to
a stationary labelling $\lambda$ on the edges of the $2-adic$ odometer
system whose image consists of the two elements $\{s,t\} \subseteq \mathfrak{S}_3$.
In this example the map $\rho_B$ of the theorem from the quotient of the dynamical system for $(A,\sigma )$ by $\mathfrak{S}_3$
to the $2-adic$ odometer system is not an isomorphism. This follows from an observation of Matui [M, see para before lemma 2.2]
and from the fact that $\mathfrak{S}_3$ is not a cyclic group.

\smallskip

\noindent{\bf 2.24. }Finally it follows from the assertion  part (II) (i) of the theorem and from [M, Lemma 2.2] that there must exist a continuous function $c:\overline{X}_{B(\lambda)}^{\mathcal Q}\longrightarrow G$ such that $(X_{B(\lambda)}^{\mathcal Q }, T_{B(\lambda)}^{\mathcal Q})
\simeq (\overline{X}_{B(\lambda)}^{\mathcal Q}\times G,\Psi )$ where 
$\Psi (y,g)=(\overline{T}_{B(\lambda)}^{\mathcal Q}(y),
gc(\overline{T}_{B(\lambda)}^{\mathcal Q}(y))$.

A function $c$ with this property is easy to describe. Let $y$ belong to level $n$ tower of
$\overline{X}_{B(\lambda)}^{\mathcal Q}$ parametrized  by the $G$-orbit of $[(u,g_1),(v,g_2),(w,g_3)]\in
\overline{V}_{n,\lambda}^{\mathcal Q}$. Define (for $n\geq 1$) $c_n(y)=g_1^{-1}g_2$ if $y$ belongs to the bottom floor of the above tower ; otherwise, define $c_n(y)=1\in G$. All the cocycles $c_n$ are cohomologous (${\it cf.}$[M, definition 2.1]).

\smallskip

\noindent{\bf 2.25. }We also remark that Matui's ({\it cf. }[M]) dynamical systems $(X\times G, \psi)$ of the form $\psi (x,g) =
(T(x), g\cdot c(T(x))$ for a $G$-valued cocycle $c$ on $X$ arise as the skew-product system of this
article associated to a particular $G$-valued labelling of the edges of a Bratteli diagram for $(X,T)$. In fact, let us use
Matui's notation in [M, p.140]. There he considers a nested sequence of K-R partitions $({\mathcal P}_n)_n$ in $X$
such that $c$ is constant on the sets of the partition and
a nested sequence of K-R partitions $({\mathcal Q}_n)_n$ in $X\times G$. They are related as follows:-
denoting, as Matui does ({\it cf. }[M,p.138]), the floors of  ${\mathcal P}_n$ in $X$ by $X(n,v,k)$,  where $X(n,v,k)$
is the $k^{th}$ floor of a level-$n$ tower ${\mathcal S}_v$,
the sets of the partition ${\mathcal Q}_n$ in $X\times G$ are
$X(n,v,k) \times \{g\}, (g\in G)$,([{\it loc.cit },p.140]).  Let $(V,E,\geq)$ be the Bratteli diagram of $({\mathcal P}_n)_n$. Let $w\in V_{n+1}$.
Let $(e_1,e_2,\cdots,e_{\ell})$ be the edges in $r^{-1}(w)$ enumerated in the linear order in $r^{-1}(w)$.
Let $h(v)$ denote the height of the K-R tower ${\mathcal S}_v$ for a vertex $v\in V$.
For $i = 1,2,\cdots,\ell$, choose $(k_i)_i$, such that
\begin{enumerate}
\item[]$1\leq k_1<k_2<k_3<\cdots<k_{\ell}  \leq h(w)$ and
\item[]$X(n+1,w,k_i) \subseteq X(n,v_i,h(v_i))$ for some $v_i\in V_n$.
\end{enumerate}
There is a unique level-$(n+1)$ tower for the K-R
partition  ${\mathcal Q}_n$ of $X\times G$ whose top floor is $X(n+1,w,h(w))\times \{e\}$ where $e\in G$ is the identity element
of $G$. Let the floors of this tower be $X(n+1,w,1)\times \{g_1\}, X(n+1,w,2)\times \{g_2\}, X(n+1,w,3)\times \{g_3\},
\cdots,X(n+1,w,h(w))\times \{g_{h(w)}\}$.  Then the labelling of $e_i$ is
$g_{k_i}$.

\smallskip

We end with a definition and a comment and indicate scope for further progress.

\noindent{\bf 2.26. Definition. } Let $\lambda$ and $\mu$ be two $G$-valued edge labellings of $(V,E,\geq)$. We say that
$\lambda$ and $\mu$ are cohomologous if there exists a function $\beta: V \to G$ such that
$$ \beta (r(e)) \mu (e) = \lambda (e) \beta (s(e)), \text{for }e\in E_n, (n\geq 1).$$
Such a map $\beta$ gives rise to an isomorphism
$$\Phi_{\beta}:(V_{\lambda}, E_{\lambda},\geq) \longrightarrow  (V_{\mu}, E_{\mu},\geq)$$
where $\Phi_{\beta}(v,g)= (v,g\beta (v)) \text{ and } \Phi_{\beta}(e,g)= (e,g\beta (r(e)))$.
This notion is related to Matui's definition as follows: recall $({\it cf. }[M, 2.1])$ that Matui calls
two $G$-valued cocycles $c,c'$ cohomologous if $\exists $  a continuous function $b: X \to G$ such that
$ c(T(x)) b(T(x))=b(x)c'(T(x))$. Let $\lambda \text{ and }\mu$ be the labellings associated to $c \text{ and }c'$, respectively
as in (2.25). We show that $\lambda$ and $\mu$ are cohomologous. Let us use the notation of (2.24). By doing a telescoping of
$(V,E,\geq)$ if necessary, we can suppose that $c,c'$ and $b$ are all constant on the sets of the K-R partition ${\mathcal P}_1$.
Then we just have to take
$$\beta (v) = b(x), \text{ where }x\in \text{ the top floor of the tower }{\mathcal S}_v$$
for a vertex $v\in V$.

\smallskip

 We do not know of any example of a substitutional system $(A_{\lambda}, \sigma_{\lambda})$ (see $2.5$) arising from a
stationary labelling $\lambda$ with values in a non-cyclic group $G$ , such that $(A,\sigma)$ and $(A_{\lambda}, \sigma_{\lambda})$
are both {\it Toeplitz flows} defined below as usual ({\it cf. eg.}[GJ,p.1695]).

\noindent{\bf 2.27. Definition. }A Toeplitz sequence is a non-periodic sequence $\eta = (\eta_n)_{n \in \Z}$ in $A^{\Z}$,
where $A$ is any finite alphabet set, so that for each $m \in \Z $ there exists $n \in \N$ so that $\eta_m = \eta_{m+kn} $ for all
$k\in \Z$.

\smallskip

 The following gives such an example with $G$ a cyclic group of order $k$. Let  $(A,\sigma)$ be a
primitive aperiodic proper substituion of constant length $n$. By [GJ, Corollary 9,p.1698] this gives rise to a 
Toeplitz flow. Assume that $n$ is congruent to 1 mod $k$. Let $z$ denote a generator for $G$. We take the
stationary labelling $\lambda$ which assigns $z^i$ to the $i^{th}$ edge ranging at any vertex of the
stationary Bratteli diagram for $(A,\sigma)$ that we described in subsection 2.5. It is not hard
to see that if $\eta = (a_i)_{i \in \Z}$ is a Toeplitz sequence for $(A,\sigma)$ then 
$\widetilde{\eta} \text{ }{\overset{def.}=}\{(a_i,z^i)\}_{i \in \Z}$ is a Toeplitz sequence for $(A_{\lambda},\sigma_{\lambda})$.

 \noindent{\bf 2.28. Further observations. }{\it After the first version of the article was submitted for publication we noticed why, for non-cyclic groups $G$, one cannot have Toeplitz sequences occurring in primitive substitutional systems 
$(A_{\lambda},\sigma_{\lambda})$ arising from $G$-valued stationary edge labellings $\lambda$ of another substitutional system
$(A,\sigma)$. }  Here's a  proof. Suppose $\widetilde{\eta} \buildrel {def.} \over {=}
(a_i,g_i)_{i \in \Z}$ is a Toeplitz sequence in  $(A_{\lambda},\sigma_{\lambda})$ (which is assumed to be primitive). Then $\forall g \in G, \widetilde{\eta}_g \buildrel {def.} \over {=} (a_i,gg_i)_{i \in \Z}$ are also
Toeplitz sequences. Further, they all admit the same period structure $p=(p_0,p_1,p_2,\cdots)$ ({\it cf.} [DKL, p.220, GJ,p.1695]). 
 More importantly, all these Toeplitz sequences lie in the same (minimal) subshift dynamical system $X_{\lambda}$ in the alphabets $A_{\lambda}$ and
the $\Z$- orbit of any one of these Toeplitz sequences $\widetilde{\eta}_g$ is dense in $X_{\lambda}$. Let $(G_p,1)$ be the
maximal uniformly continuous factor of $(X_{\lambda},T_{\lambda})$ ({\it cf.}[GJ, p.1696]) and let $\pi : (X_{\lambda},T_{\lambda}) \longrightarrow (G_p,1)$
denote the corresponding factor map. For a Toeplitz sequence $\omega, \pi^{-1}(\pi(\omega)) = \{\omega \}$.({\it cf.}[W], [GJ,p.1696],
[DL,Theorem 6,p.167]). Define elements $h_g, (g\in G)$ of the monothetic group $G_p$ by $h_g = \pi (\widetilde{\eta}_g)$. They are pairwise
distinct. For $g\in G$, denote by $\tau_g :(X_{\lambda},T_{\lambda}) \longrightarrow (X_{\lambda},T_{\lambda})$ the unique isomorphism which replaces an alphabet $(a, \theta)$ by $(a, g\theta)$, for $\theta \in G$. Thus, for example, $\tau_g (\widetilde{\eta}_{\theta})
= (\widetilde{\eta}_{g\theta})$.

 As observed in [DKL, sec.2, paragraph 1], we have the following commutative diagram
$$\begin{matrix}
(X_\lambda, T_\lambda )&\buildrel \tau_g \over \longrightarrow &(X_\lambda, T_\lambda )\\
\pi\downarrow &&\downarrow \pi \\
(G_p,1)&\buildrel {h_g } \over \longrightarrow
&(G_p,1 )
\end{matrix}$$
where the bottom map denotes group operation by $h_g$.  Composition of two such diagrams for $g_1,g_2 \in G$ yields
the diagram for $g_1g_2$.  Thus, the group $G$ sits faithfully as a finite subgroup of $G_p$. Hence, it has to be isomorphic
to a subgroup of the finite cyclic group $\Z_{p_i}$ for one of the essential periods in $p_0,p_1,p_2,\cdots,p_i,\cdots$. \hfill $\square $
 
 On the other hand, if we drop the `Toeplitz' requirement, the subshift dynamical system of 2.23
in the alphabet set  $A= \{a_{\tau} \mid \tau \in \mathfrak{S}_3 \}$ provides examples of sequences
$\eta = (a_{\tau_i} )_{i \in \Z}$ such that all the sequences $\eta_g  \buildrel {def.} \over {=} (a_{g\tau_i} )_{i \in \Z},
(\tau \in \mathfrak{S}_3)$ lie in the same (minimal) subshift dynamical system in the alphabets $A$.

%\newpage
\smallskip

\centerline{\bf References}

%\smallskip

\begin{enumerate}

\item[[DKL]] {\sc T. Downarowiz, J. Kwiatkowski and Y. Lacroix}. A criterion for Toeplitz flows to be
topologically isomorphic and applications, {\it Coll. Math.}{\bf 68}, (1995) 219-228.
\item[[DL]]{\sc T. Downarowiz and Y. Lacroix}. Almost 1-1 extensions of Furstenberg-Weiss type and 
applications to Toeplitz flows, {\it Studia Math.} {\bf 130}, (1998) 149-170.
\item[[DHS]] {\sc F. Durand, B. Host and C. Skau}. Substitutional dynamical systems, Bratteli
diagrams and dimension groups, {\it Ergodic Th. and Dynam. Sys.} {\bf 19}, (1999) 953-993.
\item[[GJ]] {\sc R. Gjerde and \O. Johansen}. Bratteli-Vershik models for Cantor minimal systems:
applications to Toeplitz flows, {\it Ergodic Th. and Dynam. Sys.} {\bf 20}, (2000) 1687-1710.
\item[[HPS]] {\sc R.H. Herman, I.F. Putnam and C.F. Skau}. Ordered Bratteli diagrams, dimension groups
and topological dynamics, {\it Internat. J. Math.} {\bf 3}, (1992) 827-864.
\item[[ KP]] {\sc A. Kumjian and D. Pask}. $C^\ast $-algebras of directed graphs and group actions,
{\it Ergodic Th. and Dynam. Sys.} {\bf 19}, (1999) 1503-1519.
\item[[  M]] {\sc H. Matui}. Finite order automorphisms and dimension groups of Cantor minimal systems, {\it J.
Math. Soc.  Japan} {\bf 54}, (2002) 135-160.
\item[[  W]] {\sc S. Williams}. Toeplitz minimal flows which are not uniquely ergodic, {\it Z. Wahrs. verw. Gebiete.} 
{\bf 67}, (1984) 95-107.

\end{enumerate}

\bigskip

\noindent {\sc Lamath,  Le Mont Houy

\noindent Universit\'e de Valenciennes

\noindent 59313 Valenciennes Cedex 9 (France)}

\noindent E-mail address : aziz.elkacimi@univ-valenciennes.fr

\bigskip

\noindent {\sc School of Mathematics

\noindent Tata Institute of Fundamental Research

\noindent Homi Bhabha Road, Colaba

\noindent 400 005 Mumbai (India)}

\noindent E-mail address : sarathy@math.tifr.res.in

\end{document}